\documentstyle[12pt]{article}

\newtheorem{teor}{Theorem}
\newtheorem{prop}{Proposition}

\newtheorem{lem}{Lemma}
\newtheorem{definition}{Definition}

\def\beq{\begin{equation}}
\def\eeq{\end{equation}}
\def\bea{\begin{eqnarray}}
\def\eea{\end{eqnarray}}
\def\beann{\begin{eqnarray*}}
\def\eeann{\end{eqnarray*}}
\def\beasn{\begin{sneqnarray}}
\def\eeasn{\end{sneqnarray}}
\def\ben{\begin{enumerate}}
\def\een{\end{enumerate}}
\def\bit{\begin{itemize}}
\def\eit{\end{itemize}}
\def\dst{\(\displaystyle}
\def\proof{( {\sl Proof} )\quad}

\def\derpar#1#2{\frac{\partial{#1}}{\partial{#2}}}

\def\mapping#1{\mathrel{\mathop{\longrightarrow}\limits^{#1}}}

\def\feble#1{\mathrel{\mathop =\limits_{#1}}}
\def\forta#1{\mathop{\mathpalette\@vereq\sim}\limits_{#1}}

\def\moment#1#2#3{{#1}_{#2}, \ldots, {#1}_{#3}}
\def\qed{\ifvmode\removelastskip\fi
{\unskip\nobreak\hfil\penalty50\hbox{}\nobreak\hfil
\hbox{\vrule height1.2ex width1.2ex}\parfillskip=0pt
\finalhyphendemerits=0 \par\smallskip}}
\def\vf{{\cal X}}
\def\df{{\mit\Omega}}
\def\Lag{{\cal L}}

\def\d{{\rm d}}

\def\Real{{\bf R}}

\def\Ker{\mathop{\rm Ker}\nolimits}

\def\inn{\mathop{i}\nolimits}

\def\Tan{{\rm T}}

\def\Lie{\mathop{\rm L}\nolimits}
\def\Cinfty{{\rm C}^\infty}

\catcode`@=12
\def\tabaddress#1{{\small\it\begin{tabular}[t]{c}#1 \\[1.2ex]\end{tabular}}}
\def\UPCMAT{\it Departamento de Matem\'atica Aplicada y Telem\'atica\\
   Edificio C-3, Campus Norte UPC\\
   C/ Jordi Girona 1\\
   E-08034 BARCELONA\\
   SPAIN}

\def\ls{((E,M;\pi),\Lag )}
\def\lag{\pounds}

\title{MULTIVECTOR FIELDS AND CONNECTIONS.
SETTING LAGRANGIAN EQUATIONS IN FIELD THEORIES}
\author{\sc Arturo Echeverr\'{\i}a-Enr\'\i quez,
   \\
\sc Miguel C. Mu\~noz-Lecanda
\thanks{{\bf e}-{\it mail}: MATMCML@MAT.UPC.ES},
   \\
\sc Narciso Rom\'an-Roy
\thanks{{\bf e}-{\it mail}: MATNRR@MAT.UPC.ES}
   \\
   \tabaddress{\UPCMAT}}
\date{ }

\begin{document}
\maketitle
\thispagestyle{empty}
\setcounter{page}{0}

\begin{abstract}
The integrability of multivector fields in a differentiable manifold is studied.
Then, given a jet bundle $J^1E\to E\to M$,
it is shown that integrable multivector fields in $E$
are equivalent to integrable connections in the bundle $E\to M$
(that is, integrable jet fields in $J^1E$).
This result is applied to the particular case of
multivector fields in the manifold $J^1E$ and connections in the bundle
$J^1E\to M$ (that is, jet fields in the repeated jet bundle $J^1J^1E$),
in order to characterize integrable multivector fields and connections whose integral
manifolds are canonical lifting of sections.

These results allow us to set the Lagrangian evolution equations
for first-order classical field theories in three
equivalent geometrical ways (in a form similar to that in which
the Lagrangian dynamical equations of non-autonomous mechanical systems
are usually given). Then, using multivector fields;
we discuss several aspects of these evolution equations
(both for the regular and singular cases);
namely: the existence and non-uniqueness of solutions,
the integrability problem and Noether's theorem;
giving insights into the differences between mechanics and field theories.
\end{abstract}


\vfill \hfill
\vbox{\raggedleft PACS: 0240, 0320, 0350. \\
AMS s.\,c.\,(1991): 53C80, 55R10, 58A20,70G50, 70H99.}\null
\clearpage

\tableofcontents

\section{Introduction}

The theory of {\sl multisymplectic manifolds}
(that is, differentiable manifolds endowed with a
closed non-degenerate $k$-form, with $k\geq 2$)
has been revealed as a powerful tool
for geometrically describing some problems in physics.
In particular, jet bundles and their duals,
equipped with specific multisymplectic forms,
are the suitable geometrical frameworks for describing
the Lagrangian and Hamiltonian formalisms
of first-order classical field theories
(a non-exhaustive list of references is
\cite{BSF-88}, \cite{CCI-91}, \cite{EM-92}, 
\cite{EMR-96}, \cite{Gc-73}, \cite{GS-73},
\cite{GIMMSY-mm}, \cite{Sd-94c}, \cite{Sa-87}).

As a consequence of this fact, the study of multisymplectic manifolds
and their properties has increased lately
\cite{CIL-96a}, \cite{CIL-96b}, \cite{IEMR-97};
in particular, those concerning the behavior of
{\sl multisymplectic Hamiltonian systems},
which are the generalization of the corresponding symplectic case.
This generalization requires the use of
{\sl multivector fields} as a fundamental tool, and their contraction with
differential forms, which is the intrinsic formulation of the
systems of partial differential equations locally describing the field.
Hence, the study of the integrability of such equations;
(that is, of the corresponding multivector fields), is of
considerable interest and constitutes the first aim in this paper.

A particular situation of special relevance arises
when multivector fields in fiber bundles are considered.
In these cases, a result on the equivalence
of integrable multivector fields and jet fields, or connections
in the corresponding first-order jet bundles, can be stated.

Furthermore, in the jet bundle description of
classical field theories, the evolution equations are usually obtained
using the multisymplectic form in order to characterize the
critical sections which are solutions of the problem.
Nevertheless, in mechanics we can write the evolution equations
in a more geometric-algebraic manner using vector fields,
and then obtain the critical sections as integral curves of these vector fields.
Different attempts have been made to achieve the same goal
with the evolution equations of field theories.
So, for the Lagrangian formalism, this is done in two different ways:
by using {\sl Ehresmann connections} in a jet bundle
\cite{LMM-95}, \cite{Sa-89} or, what is equivalent,
their associated {\sl jet fields} \cite{EMR-96}.
Moreover, in  \cite{Ka-94}, \cite{Ka-95a}, \cite{Ka-97a} and \cite{Ka-98},
the evolution equations in the Hamiltonian formalism
and generalized Poisson brackets are stated (mainly in local coordinate
terms) using multivector fields (although, as far as we know, the first use of multivector fields in the realm of field theories can be found in \cite{GS-73}).
A further goal of this work is to carry out a deeper and pure geometrical analysis
of these procedures for the Lagrangian formalism, showing that
all these ways are, in fact, equivalent.

Thus, the aims of this paper are:
\ben
\item
To study the integrability of
multivector fields in differentiable manifolds in general.
\item
To set the equivalence between jet fields
in the jet bundle $J^1E$ and a certain kind of multivector fields in $E$.
This equivalence will in fact be used for jet fields
in the ``repeated'' jet bundle $J^1J^1E$
and the corresponding multivector fields in the jet bundle $J^1E$.
\item
To use these results to state the Lagrangian evolution equations
for first-order classical field theories using multivector fields; and to
prove this formalism is equivalent to the above mentioned ways of setting these equations.
\item
To use, in particular, this {\sl multivector field formalism}
to study some features of Lagrangian field theories;
namely: existence and non-uniqueness of solutions,
the problem of the integrability and Noether's theorem.
In this way the differences between mechanics and field theories
will be made evident.
\een

The structure of the paper is the following:

The first part is devoted to an analysis of the fist two items. In particular:

In Section \ref{imvfdm} we introduce the notion and characterization of
{\sl integrable multivector fields} in a manifold, as well as
several properties for them.

Section \ref{ia} is devoted to the presentation of an algorithm for finding the
submanifold (if it exists) where a non-integrable multivector field is integrable.
This algorithm is applied to solve a
particular system of partial differential equations in Section \ref{example0}.

In Section \ref{mfjfjb} and \ref{mfcjb}, we study the integrability of
jet fields, connections and special kinds of multivector fields in jet bundles;
and their equivalence.
This equivalence is considered for the particular case of jet fields
in $J^1J^1E$, and the suitable multivector fields in $J^1E$.

In Section \ref{hmjfrjb}, after giving different geometrical characterizations for
{\sl semi-holonomic} and {\sl holonomic} multivector fields, they are finally related
with semi-holonomic and holonomic jet fields.

In the second part, the applications to Lagrangian first-order field theories
are considered. Thus:

Section \ref{eblft} is a short background on Lagrangian first-order classical field theory.

Section \ref{lfcft} is devoted to establishing
the evolution equations for classical field theories
(the {\sl Euler-Lagrange equations}) using
multivector fields in $J^1E$, showing that this is
equivalent to using jet fields in $J^1J^1E$ or
their associated connections in $J^1E\to M$.

In Sections \ref{aeerlft} and \ref{dslft}, the analysis of the evolution equations is made for regular and singular Lagrangian field theories,
using multivector fields. In particular, we show that the existence of integrable multivector fields solution of the equations is not guaranteed, (even in the regular case).
The non-uniqueness of solutions is also discussed, as well as the
integrability of these solutions.
Furthermore, algorithmic procedures for finding
multivector fields solutions in the singular case
is outlined in section \ref{dslft}.

Section \ref{ntmvf} deals with the {\sl generalized symmetries}
of a Lagrangian system, and a version of {\sl Noether's theorem}
using the multivector field formalism is proved.

An example, which is a quite general version of many typical models in field theories,
is analyzed in section \ref{example}.

Finally, a summary with the main results in the work is given.

An appendix is included as a remainder of the definitions and main properties
of the canonical structures in a jet bundle, which are used in some parts of the work.

All manifolds are real, paracompact, connected and $C^\infty$. All maps are $C^\infty$.
Sum over crossed repeated indices is understood.

\section{Multivector fields and connections}

We begin this work by studying the integrability of
multivector fields in differentiable manifolds in general;
and setting the equivalence between jet fields
and multivector fields in jet bundles.

\subsection{Multivector fields in differentiable manifolds}
\protect\label{imvfdm}

Let $E$ be a $n$-dimensional differentiable manifold.
Sections of $\Lambda^m(\Tan E)$ are called
$m$-{\sl multivector fields} in $E$.
We will denote by $\vf^m (E)$ the set of $m$-multivector fields in $E$.
In general, given a $m$-multivector field $Y\in\vf^m(E)$,
for every $p\in E$, there exists an open neighbourhood $U_p\subset E$
and $Y_1,\ldots ,Y_r\in\vf (U_p)$ such that
\beq
Y\feble{U_p}\sum_{1\leq i_1<\ldots <i_m\leq r}
f^{i_1\ldots i_m}Y_{i_1}\wedge\ldots\wedge Y_{i_m}
\label{contrac}
\eeq
with $f^{i_1\ldots i_m}\in\Cinfty (U_p)$ and $m\leq r\leq{\rm dim}\, E$.
The situation we are interested in is when,
for every $p\in E$, we can take $r=m$; therefore
\cite{CDD-amp}, \cite{BCG-91}:

\begin{definition}
A $m$-multivector field $Y\in\vf^m(E)$ is said to be
{\rm decomposable} iff there are
$Y_1,\ldots ,Y_m\in\vf (E)$ such that
$Y=Y_1\wedge\ldots\wedge Y_m$.

The multivector field $Y\in\vf^m(E)$ is said to be
{\rm locally decomposable} iff,
for every $p\in E$, there exists an open neighbourhood $U_p\subset E$
and $Y_1,\ldots ,Y_m\in\vf (U_p)$ such that
$Y\feble{U_p}Y_1\wedge\ldots\wedge Y_m$.
\end{definition}

Every multivector field $Y\in\vf^m(E)$ defines a derivation $\inn (Y)$
of degree $-m$ in the algebra of differential forms $\df (E)$.
If $\Omega\in\df^k(E)$ is a differentiable $k$-form in $E$;
using (\ref{contrac}), this is
$$
\inn(Y)\Omega\feble{U_p}
\sum_{1\leq i_1<\ldots <i_m\leq r}f^{i_1\ldots i_m}
\inn(Y_1\wedge\ldots\wedge Y_m)\Omega =
\sum_{1\leq i_1<\ldots <i_m\leq r}f^{i_1\ldots i_m}
\inn (Y_1)\ldots\inn (Y_m)\Omega
$$
if $k\geq m$, and it is obviously equal to zero if $k<m$.
The $k$-form $\Omega$ is said to be {\sl $j$-nondegenerate}
(for $1\leq j\leq k-1$) iff, for every $p\in E$ and $Y\in\vf^j(E)$,
$\inn(Y_p)\Omega_p =0\ \Leftrightarrow \ Y_p=0$
(see \cite{CIL-96a}, \cite{CIL-96b} for more information on these topics).

Let $D$ be a $m$-dimensional distribution in $E$;
that is, a $m$-dimensional subbundle of $\Tan E$.
Obviously sections of $\Lambda^mD$ are $m$-multivector fields in $E$.
The existence of a non-vanishing global section of $\Lambda^mD$
is equivalent to the orientability of $D$.
We are interested in the relation between non-vanishing
$m$-multivector fields in $E$ and $m$-dimensional distributions
in $\Tan E$.
So we set the following:

\begin{definition}
A non-vanishing $m$-multivector field $Y\in\vf^m(E)$ and
a $m$-dimensional distribution $D\subset\Tan E$
are {\rm locally associated} iff there exists a connected open set
$U\subseteq E$ such that $Y\vert_U$ is a section of $\Lambda^mD\vert_U$.
\end{definition}

As a consequence of this definition,
if $Y,Y'\in\vf^m(E)$ are non-vanishing multivector fields
locally associated, on the same connected open set $U$,
with the same distribution $D$, then there exists a
non-vanishing function $f\in\Cinfty (U)$ such that
$Y'\feble{U}fY$. This fact defines an equivalence relation in the
set of non-vanishing $m$-multivector fields in $E$, whose equivalence classes
will be denoted by $\{ Y\}_U$. Therefore we can state:

\begin{teor}
There is a bijective correspondence between the set of $m$-dimensional
orientable distributions $D$ in $\Tan E$ and the set of the
equivalence classes $\{ Y\}_E$ of non-vanishing, locally decomposable
$m$-multivector fields in $E$.
\label{bijcor}
\end{teor}
\proof
Let $\omega\in\df^m(E)$ be an orientation form for $D$.
If $p\in E$ there exists an open neighbourhood $U_p\subset E$
and $Y_1,\ldots ,Y_m\in\vf (U_p)$,
with $\inn(Y_1\wedge\ldots\wedge Y_m)\omega >0$, such that
\dst D\vert_{U_p}={\rm span}\, \{Y_1,\ldots ,Y_m\}\) .
Then $Y_1\wedge\ldots\wedge Y_m$ is a representative of a class
of $m$-multivector fields associated with $D$ in $U_p$.
But the family $\{ U_p\ ;\ p\in E\}$ is a covering of $E$;
let $\{ U_\alpha\ ;\ \alpha\in A\}$ be
a locally finite refinement and $\{ \rho_\alpha\ ;\ \alpha\in A\}$
a subordinate partition of unity.
If $Y^\alpha_1,\ldots ,Y^\alpha_m$ is a local basis of $D$ in $U_\alpha$,
with $\inn(Y^\alpha_1\wedge\ldots\wedge Y^\alpha_m)\omega >0$,
then \dst Y=\sum_\alpha\rho_\alpha Y^\alpha_1\wedge\ldots\wedge Y^\alpha_m\)
is a global representative of the class of non-vanishing
$m$-multivector fields associated with $D$ in $E$.

The converse is trivial because if
$Y\vert_U=Y^1_1\wedge\ldots\wedge Y^1_m=Y^2_1\wedge\ldots\wedge Y^2_m$,
for different sets $\{ Y^1_1,\ldots ,Y^1_m\}$, $\{Y^2_1,\ldots ,Y^2_m\}$,
then
$span\, \{ Y^1_1,\ldots ,Y^1_m\} =span\, \{ Y^2_1,\ldots ,Y^2_m\}$.
\qed

{\bf Comments}:
\bit
\item
If $Y\in\vf^m(E)$ is a non-vanishing $m$-multivector field
and $U\subseteq E$ is a connected open set, the distribution associated with
the class $\{ Y\}_U$ will be denoted by ${\cal D}_U(Y)$.
If $U=E$ we will write simply ${\cal D}(Y)$.
\item
If $D$ is a non-orientable $m$-dimensional distribution in $\Tan E$,
for every $p\in E$, there exist an open neighbourhood $U_p\subset E$
and a non-vanishing multivector field $Y\in\vf^m(U)$ such that
$D\vert_{U_p}={\cal D}_U(Y)$.
\eit

\begin{definition}
Let $Y\in\vf^m(E)$ a multivector field.
\ben
\item
A submanifold $S\hookrightarrow E$, with ${\rm dim}\, S=m$,
is said to be an {\rm integral manifold} of $Y$ iff,
for every $p\in S$, $Y_p$ spans $\Lambda^m\Tan_pS$.
\item
$Y$ is said to be an {\rm integrable multivector field}
on an open set $U\subseteq E$ iff,
for every $p\in U$, there exists an
integral manifold $S\hookrightarrow U$ of $Y$, with $p\in S$.

$Y$ is said to be {\rm integrable} iff it is integrable in $E$.
\een
\label{defint}
\end{definition}

Obviously, every integrable multivector field is non-vanishing.

Now, bearing in mind the statement in theorem \ref{bijcor},
we can define:

\begin{definition}
Let $Y\in\vf^m(E)$ be a multivector field.
\ben
\item
$Y$ is said to be {\rm involutive} on a connected open set $U\subseteq E$ iff
it is locally decomposable in $U$ and its associated distribution
${\cal D}_U(Y)$ is involutive.
\item
$Y$ is said to be {\rm involutive} iff it is involutive on $E$.
\item
$Y$ is said to be {\rm locally involutive} around $p\in E$ iff
there is a connected open neighbourhood $U_p\ni p$
such that $Y$ is involutive on $U_p$.
\een
\end{definition}

Note that if $Y$ is locally involutive around every $p\in E$,
then it is involutive.

Now, the classical {\sl Frobenius' theorem} can be reformulated
in this context as follows:

\begin{prop}
A non-vanishing and locally decomposable multivector field
$Y\in\vf^m(E)$ is integrable on a
connected open set $U\subseteq E$ if, and only if, it is involutive on $U$. 
\label{integr}
\end{prop}

Note that if a multivector field $Y$ is integrable,
then so is every other in its equivalence class $\{ Y\}$,
and all of them have the same integral manifolds.

As is well known, a $m$-dimensional distribution $D$ is integrable if,
and only if, it is locally spanned by a set of vector fields
$(Y_1,\ldots ,Y_m)$ such that $[Y_\mu,Y_\nu]=0$, for every pair $Y_\mu,Y_\nu$.
Thus a multivector field $Y\in\vf^m(E)$ is integrable if, and only if,
for every $p\in E$, there exists an open neighbourhood $U_p\subset E$
and $Y_1,\ldots ,Y_m\in\vf (U_p)$ such that:
\ben
\item
$Y_1,\ldots ,Y_m$ span ${\cal D}_{U_p}(Y)$.
\item
$[Y_\mu,Y_\nu]=0$, for every pair $Y_\mu,Y_\nu$.
\een
Then there is a non-vanishing function $f\in\Cinfty (U_p)$ such that
\dst Y\vert_{U_p}=fY_1\wedge\ldots\wedge Y_m\) .

Now, let $\{\tau^\mu\}$ be the local one-parameter groups of diffeomorphisms
of $Y_\mu$ around $p$ (for $\mu=1,\ldots ,m$). The second condition above
implies that $\tau^\mu\circ\tau^\nu=\tau^\nu\circ\tau^\mu$,
for every $\mu,\nu$. Then, if $S_p$ is the integral manifold of $Y$
through $p$, there exist open neighbourhoods
$V(0)\subset\Real^m$ and $W_p\subset U_p\subset E$, and a map
$$
\begin{array}{cccccc}
\tau&\colon&V(0)\times W_p&
\begin{picture}(30,5)(0,0)
\put(0,4){\vector(1,0){30}}
\end{picture}
& U_p &
\\
& &(t,q)&\mapsto&\tau_t(q):=(\tau^1_{t_1}\circ\ldots\circ\tau^m_{t_m})(q)
&\qquad (t=(t_1,\ldots ,t_m))
\end{array}
$$
verifying that:
\ben
\item
$\tau$ is injective.
\item
$\tau_{t+s}=\tau_t\circ\tau_s$, for $t,s,t+s\in V(0)$.
\item
For every $q\in W_p$, the set $\{\tau_t(q)\ ;\ t\in V(0)\}$
is an open neighbourhood of $q$ in $S_q$.
\een
The map $\tau$ is called the {\sl $m$-flow} associated with the
multivector field $Y$ (see \cite{Co-93} for the terminology
and notation). Therefore we can define:

\begin{definition}
A multivector field $Y\in\vf^m(E)$ is said to be a
{\rm dynamical multivector field} iff 
\ben
\item
$Y$ is integrable.
\item
For every $p\in E$, there exists an open neighbourhood $U_p\subset E$
and $Y_1,\ldots ,Y_m\in\vf (U_p)$ such that
$[Y_\mu,Y_\nu]=0$, for every pair $Y_\mu,Y_\nu$, and
\dst Y\vert_{U_p}=Y_1\wedge\ldots\wedge Y_m\) .
\een
\label{dynmvf}
\end{definition}

Thus, as a consequence of all these comments we have that:

\begin{prop}
Let $\{ Y\}$ be a class of integrable $m$-multivector fields.
Then there is a representative $Y$ of the class which is
a dynamical multivector field.
\end{prop}

\subsection{Integrability algorithm for locally decomposable multivector fields}
\protect\label{ia}

Every locally decomposable multivector field $Y\in\vf^m(E)$ defines
a system of (first-order) partial differential equations,
whose solutions (if they exist) give locally the integral manifolds of $Y$.

Nevertheless, in many applications, we have a locally decomposable
multivector field which is not integrable in $E$, but integrable
in a submanifold of $E$. What this means is that the corresponding
system of partial differential equations has no solution everywhere in $E$,
but only on the points of that submanifold.
For instance, this problem arises when we look for solutions
of evolution equations in field theories  (as we will see later) .

Next we present an algorithm which allows us to find this submanifold.
Although we will maintain the global notation ($E$ for the manifold),
this is a local algorithm; that is, we are in fact working on suitable open sets
in $E$. Hence, let $Y\equiv\bigwedge_{\mu=1}^m Y_\mu$.
\bit
\item
{\sl Integrability condition}:

The condition for $Y$ to be integrable is equivalent to demanding that
the distribution spanned by $\moment{Y}{1}{m}$ is involutive.
Then, let $\moment{Z}{1}{n-m}\in\vf (E)$, such that
$\moment{Y}{1}{m}, \moment{Z}{1}{n-m}$ are a local basis of $\vf (E)$.
Therefore, for every couple $Y_\mu,Y_\nu$ ($1\leq\mu,\nu\leq m$) we have
$$
[Y_\mu,Y_\nu]=\xi_{\mu\nu}^\rho Y_\rho+\zeta_{\mu\nu}^lZ_l
$$
for some functions $\xi_{\mu\nu}^\rho,\zeta_{\mu\nu}^l$.
Consider the system $\zeta^l_{\mu\nu}=0$ and let
$$
E_1=\{ p\in E \ ;\  \zeta_{\mu\nu}^l(p)=0\ ,\ \forall\mu,\nu,l\}
$$
we have three options:
\ben
\item
$E_1=E$.

Then the distribution spanned by $\moment{Y}{1}{m}$ is involutive,
and the multivector field $Y$ is integrable in $E$.
If, in addition, all the functions $\xi_{\mu\nu}^\rho$ are zero everywhere,
then $Y$ is also a dynamical multivector field.
\item
$E_1=\emptyset$.

Then the distribution spanned by $\moment{Y}{1}{m}$ is not involutive at any point
in $E$, and hence the multivector field $Y$ is not integrable.
\item
$E_1$ is a proper subset of $E$.

In this case we assume that $E_1$ is a closed submanifold of $E$
and the functions $\zeta_{\mu\nu}^l$ are the constraints locally defining $E_1$.
The distribution spanned by $\moment{Y}{1}{m}$ is involutive on $E_1$;
that is, the multivector field $Y$ is involutive on $E_1$.
Moreover, if all the functions $\xi_{\mu\nu}^\rho$ are zero everywhere in $E_1$,
then $Y$ is also a dynamical multivector field in $E_1$.

Obviously, $Y\colon E_1\to\Lambda^m\Tan E\vert_{E_1}$.
If, in addition, $Y\colon E_1\to\Lambda^m\Tan E_1$,
then $Y$ is an integrable multivector field tangent to $E_1$.
Nevertheless, this is not the case in general,
so we need the following:
\een
\item
{\sl Tangency condition}:

Consider the set
$$
E_2:=\{ p\in E_1 \ ;\ Y(p)\in\Lambda^m\Tan_{p}E_1\}
$$
For $E_2$ we have the same problem, so we define inductively, for $i>1$,
$$
E_i:=\{ p\in E_{i-1} \ ;\ Y(p)\in\Lambda^m\Tan_{p}E_{i-1}\}
$$
and assume that we obtain a sequence of non-empty (closed) submanifolds
$\ldots \subset E_i\subset\ldots \subset E_1\subset E$.

Observe that the necessary and sufficient condition
for a locally decomposable multivector field
$Y=Y_1\wedge\ldots\wedge Y_m$ to be tangent to $E_i$
is that $Y_\mu$ is tangent to $E_i$, for every $\mu$.
Hence, if \dst E_f:=\cap_{1\geq 1}E_i\) , and assuming that $E_f$ is a non-empty (closed)
submanifold of $E$ with $\dim\, E_f\geq m$, we have that
$Y\colon E_f\to\Lambda^m\Tan E_f$;
this is involutive and therefore integrable on the manifold $E_f$.

Using the constraints, we have that, if $\{\zeta_{\alpha_i}^{(i)}\}$ is a basis of constraints
defining locally $E_i$ in $E_{i-1}$, the tangency condition is
$$
0 \feble{E_i} Y_\mu(\zeta^{(i)}_{\alpha_i}) \quad ,\quad \forall \mu,\alpha_i
$$
that is, we have
$$
E_{i+1}:=\{ p\in E_i \ ;\ Y_\mu(\zeta^{(i)}_{\alpha_i})(p)=0 \ ,\ \forall \mu,\alpha_i\}
$$
and again we have the same situation as in the case of $E_1\subset E$.
\eit

This algorithm ends in one of the following two options:
\ben
\item
We obtain a submanifold $E_f\hookrightarrow E$, with $\dim\, E_f\geq m$,
where $Y$ is an integrable multivector field (and, maybe, dynamical).
\item
We obtain a submanifold $E_f$ with $\dim\, E_f<m$, or the empty set.
Then the problem has no solution in $E$.
\een

We will call this procedure the {\sl integrability algorithm} for multivector fields.

{\bf Remark}:
As is clear, although we have developed this algorithm for
locally decomposable multivector fields, the procedure is obviously
applicable for distributions in general.

\subsection{An example}
\protect\label{example0}

In order to show how the above algorithm of integrability is applied
we take the following example of a system of partial differential equations
which does not verify the integrability condition (see \cite{Al-74})
\beann
\derpar{y_1}{x_1}=y_1-x_1-x_2 \quad & , & \quad
\derpar{y_1}{x_2}={y_1}^2-{x_1}^2-{x_2}^2-2x_1-2x_2-2x_1x_2
\\
\derpar{y_2}{x_1}=-y_2+x_1+x_2 \quad & , & \quad
\derpar{y_2}{x_2}=1
\eeann
For this system we have $E\equiv \Real^4$, with local coordinates $(x_1,x_2,y_1,y_2)$.
The associated distribution is locally spanned by the vector fields
\beann
Y_1 &=& \derpar{}{x_1}+(y_1-x_1-x_2)\derpar{}{y_1}+(-y_2+x_1+x_2)\derpar{}{y_2}
\\
Y_2 &=& \derpar{}{x_2}+({y_1}^2-{x_1}^2-{x_2}^2-2x_1-2x_2-2x_1x_2)\derpar{}{y_1}+\derpar{}{y_2}
\eeann
and a representative of the corresponding class of
locally decomposable multivector fields $\{ Y\}\subset\vf^2(E)$ is $Y=Y_1\wedge Y_2$.

\bit
\item
{\sl Integrability condition}:

As local basis of $\vf (E)$ we take the set
\dst\left(Y_1,Y_2,\derpar{}{y_1},\derpar{}{y_2}\right)\) .
Then the condition of involutivity of the above distribution leads to
$$
[Y_1,Y_2]=((x_1+x_2-y_1)^2-1)\derpar{}{y_1}=0
$$
which gives the constraint
$\zeta^{(1)}\equiv (x_1+x_2-y_1)^2-1=(x_1+x_2-y_1-1)(x_1+x_2-y_1+1)
\equiv\zeta^{(1)}_a\zeta^{(1)}_b$;
and we have the two disjoint submanifolds
\beann
E_1^a:=\{ (x_1,x_2,y_1,y_2)\ \vert \ x_1+x_2-y_1-1=0\}
\\
E_1^b:=\{ (x_1,x_2,y_1,y_2)\ \vert \ x_1+x_2-y_1+1=0\}
\eeann
and the tangency condition must hold separately for each one of them.
\item
{\sl Tangency condition}:
\ben
\item
{\sl Tangency condition for $E_1^a$}:
We obtain that $Y_1(\zeta^{(1)}_a)\feble{E_1^a}2$;
then $Y_1$ is not tangent to $E_1^a$ at any point, and
the system has not any solution on $E_1^a$.
\item
{\sl Tangency condition for $E_1^b$}:
In this case we have that
$Y_1(\zeta^{(1)}_a)\feble{E_1^b}0$ and $Y_2(\zeta^{(1)}_a)\feble{E_ 1^b}0$.
Then $Y\vert_{E_1^b}\in\vf^2(E_1^b)$ and it is integrable on $E_1^b$.
\een
\item
{\sl Integration of the system on $E_1^b$}:

The general solution of the system on $E_1^b$ is
$$
y_1=x_1+x_2+1 \quad , \quad
y_2=x_1+x_2-1+ce^{-x_1} \qquad (c\in\Real)
$$
and there exists an integral manifold through every point on $E_1^b$.
(Compare with \cite{Al-74}, pp. 47,53).
\eit

\subsection{Multivector fields in fiber bundles}
\protect\label{mfjfjb}

The particular situation in which we are interested
is the study of multivector fields in fiber bundles and,
in particular, in jet bundles.
Bearing this in mind, first we state some results concerning
multivector fields in fiber bundles in general.

Let $\pi\colon E\to M$ be a fiber bundle with $\dim\, M=m$.
We denote by ${\rm V}(\pi)$ the vertical bundle associated with $\pi$
(that is, ${\rm V}(\pi)=\Ker\Tan\pi$), and by
$\vf^{{\rm V}(\pi)}(E)$ the corresponding sections or vertical vector fields.

We are interested in the case that the integral manifolds of multivector fields
are sections of the projection $\pi$.
In order to characterize this kind of multivector fields we set:

\begin{definition}
A multivector field $Y\in\vf^m(E)$ is said to be
{\rm transverse to the projection} $\pi$ (or {\rm $\pi$-transverse})
iff, at every point $y\in E$,
$(\inn (Y)(\pi^*\omega))_y\not= 0$, for every $\omega\in\df^m(M)$
with $\omega (\pi(y))\not= 0$.
\end{definition}

{\bf Comments}:
\bit
\item
This condition is equivalent to $(\Lambda^m\Tan\pi\circ Y)(y)\not= 0$,
for every $y\in E$.
\item
If $Y\in\vf^m(E)$ is locally decomposable,
then $Y$ is $\pi$-transverse if, and only if,
$\Tan_y\pi({\cal D}(Y))=\Tan_{\pi (y)}M$, for every $y\in E$.
(Remember that ${\cal D}(Y)$ is the $m$-distribution associated to $Y$).
\eit

Then we can state:

\begin{teor}
Let $Y\in\vf^m(E)$ be integrable.
Then $Y$ is $\pi$-transverse if, and only if,
its integral manifolds are local sections of $\pi\colon E\to M$.
\label{insecmvf1}
\end{teor}
\proof
Consider $y\in E$, with $\pi(y)=x$. In a neighbourhood of $y$ there exist
$Y_1,\ldots ,Y_m\in\vf (E)$ such that $Y_1,\ldots ,Y_m$ span ${\cal D}(Y)$ and
$Y=Y_1\wedge\ldots\wedge Y_m$. But, as $Y$ is $\pi$-transverse,
$(\inn (Y)(\pi^*\omega))_y\not= 0$,
for every $\omega\in\df^m(M)$ with $\omega (x)\not= 0$.
Thus, taking into account the second comment above,
${\cal D}(Y)$ is a $\pi$-transverse distribution and
$Y_\mu\not\in\vf^{{\rm V}(\pi)}(E)$ at any point, for every $Y_\mu$.
Now, let $S\hookrightarrow E$ be the integral manifold of ${\cal D}(Y)$
passing through $y$, then
$\Tan_yS={\rm span}\{ (Y_1)_y,\ldots ,(Y_m)_y\}$.
As a consequence of all of this, and again taking into account
the second comment above, for every point $y\in S$,
$\Tan_y\pi({\cal D}(Y))=\Tan_{\pi (y)}M$,
then $\pi\vert_S$ is a local diffeomorphism
and $S$ is a local section of $\pi$.

The converse is obvious.
\qed

Observe that, in this case, if $\phi\colon U\subset M\to E$
is a local section with $\phi (x)=y$ and $\phi (U)$ is
the integral manifold of $Y$ through $y$,
then $\Tan_y({\rm Im}\,\phi)$ is ${\cal D}_y(Y)$.

Now, considering the diagram
$$
\begin{array}{ccccccc}
& & & \Lambda^m\Tan\pi & & &
\\
& \Lambda^m\Tan (\pi^{-1}(U))& &
\begin{picture}(50,10)(0,0)
\put(0,8){\vector(1,0){50}}
\put(50,2){\vector(-1,0){50}}
\end{picture}
& &\Lambda^m\Tan U&
\\
& & & \Lambda^m\Tan\phi & & &
\\
\begin{picture}(25,40)
\put(0,17){\mbox{$\Lambda^m\tau_E$}}
\end{picture}
&
\begin{picture}(10,40)(0,0)
\put(2,40){\vector(0,-1){40}}
\put(8,0){\vector(0,1){40}}
\end{picture}
&
\begin{picture}(5,40)
\put(0,17){\mbox{$Y$}}
\end{picture}
& & &
\begin{picture}(10,40)(0,0)
\put(5,40){\vector(0,-1){40}}
\end{picture}
&
\begin{picture}(15,40)
\put(0,17){\mbox{$\Lambda^m\tau_U$}}
\end{picture}
\\
& & & \phi & & &
\\
& \pi^{-1}(U) & &
\begin{picture}(50,10)(0,0)
\put(50,8){\vector(-1,0){50}}
\put(0,2){\vector(1,0){50}}
\end{picture}
& & U\subset M &
\\
& & & \pi & & &
\end{array}
$$
we have that:

\begin{prop}
$Y\in\vf^m(E)$ is integrable and $\pi$-transverse if, and only if, for every point $y\in E$,
there exists a local section $\phi\colon U\subset M\to E$
such that $\phi(\pi(y))=y$;  and a non-vanishing
function $f\in\Cinfty (E)$ such that
$\Lambda^m\Tan\phi =fY\circ\phi\circ\Lambda^m\tau_U$.
\label{transint}
\end{prop}
\proof
If $Y$ is integrable and $\pi$-transverse, then by theorem \ref{insecmvf1},
for every $y\in E$, with $\pi(y)=x$,
there is an integral local section $\phi\colon U\subset M\to E$
of $Y$ at $y$ such that $\phi (x)=y$.
Then, by definition \ref{defint}, $Y_y$ spans $\Lambda^m\Tan_y\phi$,
and hence the relation in the statement holds and the above
diagram becomes commutative on the open set $U$.

Conversely, if the relation holds, then
${\rm Im}\,\phi$ is an integral manifold of $Y$ at $y$,
then $Y$ is integrable and, as $\phi$ is a section of $\pi$,
$Y$ is necessarily $\pi$-transverse.
\qed

\subsection{Multivector fields and connections in jet bundles}
\protect\label{mfcjb}

Next we are going to establish the relation between multivector fields
and connections in jet bundles.
Let $\pi\colon E\to M$ be a fiber bundle,
$\pi^1\colon J^1E\to E$ the corresponding first-order jet bundle
of local sections of $\pi$,
and $\bar\pi^1=\pi \circ \pi^1\colon J^1E \longrightarrow M$.
If $(x^\mu,y^A)$ is a local system of coordinates adapted to the projection $\pi$
($\mu = 1,\ldots,m$, $A= 1,\ldots,N$), let $(x^\mu ,y^A,v^A_\mu)$ be the
induced local system of coordinates in $J^1E$.

With the aim of relating the elements in $\vf^m(E)$ to connections
in $\pi\colon E\to M$, we briefly recall several ways of giving a connection
(see \cite{Sa-89}):

\begin{definition}
A {\rm connection} in $\pi\colon E\to M$ is one of the following equivalent elements:
\ben
\item
A global section of $\pi^1\colon J^1E\to E$,
(that is, a mapping $\Psi\colon E\to J^1E$
such that $\pi^1\circ\Psi ={\rm Id}_E$).
It is called a {\rm jet field}.
\item
A subbundle ${\rm H}(E)$ of $\Tan E$ such that
$\Tan E={\rm V}(\pi )\oplus{\rm H}(E)$.
It is called a {\rm horizontal subbundle}.
\item
A $\pi$-semibasic $1$-form $\nabla$ on $E$ with values in $\Tan E$,
such that $\nabla^*\alpha =\alpha$, for every
$\pi$-semibasic form $\alpha\in\df^1(E)$.
It is called the {\rm connection form} or {\rm Ehresmann connection}.
\een
\end{definition}

The equivalence among these elements can be seen in \cite{Sa-89},
or it suffices to observe their local expressions.
In a natural chart $(x^\mu ,y^A,v^A_\mu)$ we have
\bea
\Psi &=& (x^\mu ,y^A,{\mit\Gamma}_\rho^A(x^\mu ,y^A))
\nonumber \\
\nabla &=&
\d x^\mu\otimes\left(\derpar{}{x^\mu}+
{\mit\Gamma}_\mu^A(x^\mu ,y^A)\derpar{}{y^A}\right)
\nonumber \\
{\rm H}(E) &=&
{\rm span}\,\left\{\derpar{}{x^\mu}+{\mit\Gamma}_\mu^A(x^\mu ,y^A)\derpar{}{y^A}\right\}
\label{locexp0}
\eea

{\bf Comments}:
\bit
\item
The horizontal subbundle ${\rm H}(E)$ is also denoted ${\cal D}(\Psi)$ when
it is considered as the distribution associated with $\Psi$.
\item
A jet field $\Psi\colon E\to J^1E$
(resp. an Ehresmann connection $\nabla$) is said to be {\rm orientable}
iff ${\cal D}(\Psi)$ is an orientable distribution on $E$.
\item
If $M$ is orientable, then every connection is also orientable.
\eit

Let $\Psi\colon E\to J^1E$ be a jet field.
A section $\phi\colon M\to E$ is said to be an
{\sl integral section} of $\Psi$ iff
$\Psi\circ\phi =j^1\phi$
(where $j^1\phi\colon M\to J^1E$ denotes the canonical lifting of $\phi$).
$\Psi$ is said to be an {\sl integrable jet field} iff it admits
integral sections through every point of $E$ or,
what is equivalent, if, and only if,
${\cal D}(\Psi )$ is an involutive distribution
(that is, ${\cal D}(\Psi )$ is integrable).
Locally, if $\phi =(x^\mu ,f^A(x^\nu ))$,
then $\phi$ is an integral section of $\Psi$ if, and only if,
$\phi$ is a solution of the following system of
partial differential equations
\beq
\derpar{f^A}{x^\mu}={\mit\Gamma}_\mu^A\circ\phi
\label{condin}
\eeq

As is well known, a jet field $\Psi$ is integrable if, and only if,
the curvature of the connection form $\nabla$ associated with $\Psi$  vanishes identically;
where the {\sl curvature} of a connection $\nabla$ is the (2,1)-tensor field defined by
$$
{\cal R}(Z_1,Z_2):=
({\rm Id}-\nabla )([\nabla (Z_1),\nabla (Z_2)])=
\inn ([\nabla (Z_1),\nabla (Z_2)])({\rm Id}-\nabla )
$$
for every $Z_1,Z_2\in\vf (E)$.
Using the coordinate expressions of the
connection form $\nabla$ and the 1-jet field $\Psi$,
a simple calculation leads to
$$
{\cal R} = \frac{1}{2}
\left(\derpar{{\mit\Gamma}_\eta^B}{x^\mu}-
\derpar{{\mit\Gamma}_\mu^B}{x^\eta}+
{\mit\Gamma}_\mu^A\derpar{{\mit\Gamma}_\eta^B}{y^A}-
{\mit\Gamma}_\eta^A\derpar{{\mit\Gamma}_\mu^B}{y^A}\right)
(\d x^\mu\wedge\d x^\eta )\otimes\derpar{}{y^B}
$$
Hence, from this local expression of ${\cal R}$ we obtain the local
integrability conditions of the equations (\ref{condin}).

Now, the relation between multivector fields and connections (jet fields)
is given as a particular case of theorem \ref{bijcor}:

\begin{teor}
There is a bijective correspondence between the set of
orientable jet fields $\Psi\colon E\to J^1E$
(that is, the set of orientable connections $\nabla$ in $\pi\colon E\to M$)
and the set of the equivalence classes
of locally decomposable and $\pi$-transverse multivector fields
$\{ Y\}\subset\vf^m(E)$.
They are characterized by the fact that
${\cal D}(\Psi)={\cal D}(Y)$.

In addition, the orientable jet field $\Psi$ is integrable if, and only if, so is $Y$,
for every $Y\in\{ Y\}$.
\label{vfmvf1}
\end{teor}
\proof
If $\Psi$ is an orientable jet field in $J^1E$,
let ${\cal D}(\Psi)$ its horizontal distribution.
Then, taking $D\equiv{\cal D}(\Psi)$, we construct
$\{ Y\}$ by applying theorem \ref{bijcor}
and, since the distribution ${\cal D}(\Psi)$ is $\pi$-transverse,
the result follows immediately.
The proof of the converse statement is similar.

Moreover, $\Psi$ is integrable if, and only if,
${\cal D}(\Psi)={\cal D}(Y)$ also is. Therefore it follows that $Y$ is also integrable,
for $Y\in\{ Y\}$,
and conversely.
\qed

As is obvious, recalling the local expression
(\ref{locexp0}), we obtain the following local expression for
a particular representative multivector field $Y$ of the class
$\{ Y\}$ associated with the jet field $\Psi$
$$
Y\equiv\bigwedge_{\mu=1}^m Y_\mu=\bigwedge_{\mu=1}^m
\left(\derpar{}{x^\mu}+{\mit\Gamma}_\mu^A\derpar{}{y^A}\right)
$$
and $\phi =(x^\mu ,f^A (x^\nu ))$
is an integral section of $Y$ if, and only if,
$\phi$ is a solution of the system of partial differential equations
(\ref{condin}).

{\bf Remarks}:
\bit
\item
If $M$ is an orientable manifold and
$\omega\in\df^m(M)$ is the volume form in $M$,
then the above representative can be obtained from the relation
\dst\inn (Y)(\pi^*\omega )=
\inn\left(\derpar{}{x^1}\wedge\ldots\wedge\derpar{}{x^m}\right)(\pi^*\omega )\)
\item
It is interesting to point out that,
if we consider the integrability algorithm described in section \ref{ia}
for this kind of multivector fields, then the integrability condition
(given by the involutivity condition for the associated distribution)
is equivalent now to the condition ${\cal R}=0$
(${\cal R}$ being the curvature of the associated jet field).
Furthermore, in this case, at the end of the algorithm, we will obtain directly
a {\sl dynamical multivector field} as the representative of the class
of integrable multivector fields.
\eit

In Lagrangian classical field theory we are interested in multivector fields
in the bundle $\bar\pi^1\colon J^1E\to M$.
Hence, all the above considerations must be adapted to the
bundle projections $J^1J^1E\mapping{\pi^1_1}J^1E\mapping{\bar\pi^1}M$.
In this way we have:
\ben
\item
Let $X\in\vf^m(J^1E)$ be integrable.
Then $X$ is $\bar\pi^1$-transverse if, and only if,
its integral manifolds are local sections of $\bar\pi^1\colon J^1E\to M$,
(theorem \ref{insecmvf1}).
\item
$X\in\vf^m(J^1E)$ is integrable and $\bar\pi^1$-transverse if, and only if, for every point
$\bar y\in J^1E$, there exists a local section $\psi\colon U\subset M\to J^1E$
such that $\psi(\bar\pi^1(\bar y))=\bar y$;
and a non-vanishing function $f\in\Cinfty (J^1E)$ such that
$\Lambda^m\Tan\psi =fX\circ\psi\circ\Lambda^m\tau_U$,
(proposition \ref{transint}).
\item
In order to relate multivector fields in $J^1E$ with connections in $\bar\pi^1\colon J^1E\to M$,
let ${\cal Y}\colon J^1E\to J^1J^1E$ be a jet field,
and $\nabla$ and ${\rm H}(J^1E)$ its associated connection form and horizontal subbundle respectively. If  $(x^\mu ,y^A,v_\mu^A)$ is a natural local chart in $J^1E$,
then the induced chart in $J^1J^1E$ is denoted
$(x^\mu ,y^A,v_\mu^A,a_\rho^A,b_{\rho\mu}^A)$.
We have the following local expressions for these elements
\beann
{\cal Y}&=&(x^\mu ,y^A,v_\mu^A,F_\rho^A(x,y,v),G_{\mu\rho}^A(x,y,v))
\\
\nabla&=&\d x^\mu\otimes
\left(\derpar{}{x^\mu}+F_\mu^A\derpar{}{y^A}+
G_{\mu\rho}^A\derpar{}{v_\rho^A}\right)
\\
{\rm H}(J^1E)&=&
{\rm span}\,\left\{\derpar{}{x^\mu}+F_\mu^A\derpar{}{y^A}+
G_{\mu\rho}^A\derpar{}{v_\rho^A}\right\}
\eeann
\item
Let $\psi =(x^\mu ,f^A (x),g^A_\mu (x))$ be a local section of $\bar\pi^1$.
It is an integral section of ${\cal Y}$ if, and only if,
$\psi $ is a solution of the following system of partial differential equations
\beq
\derpar{f^A}{x^\mu}=F_\mu^A\circ\psi
\qquad \derpar{g_\rho^A}{x^\mu}=G_{\rho\mu}^A\circ\psi
\label{sisteq}
\eeq
and the necessary and sufficient condition for this system to be integrable is that
the curvature ${\cal R}$ of ${\cal Y}$ vanishes everywhere;
that is, in coordinates
\bea
0 =
\left(\derpar{F_\eta^B}{x^\mu}+F_\mu^A\derpar{F_\eta^B}{y^A}+
G_{\gamma\mu}^A\derpar{F_\eta^B}{v_\gamma^A}-
\derpar{F_\mu^B}{x^\eta}-F_\eta^A\derpar{F_\mu^B}{y^A}-
G_{\rho\eta}^A\derpar{F_\mu^B}{v_\rho^A}\right)
(\d x^\mu\wedge\d x^\eta )\otimes\derpar{}{y^B}+
\nonumber
\\
\left(\derpar{G_{\rho\eta}^B}{x^\mu}+F_\mu^A\derpar{G_{\rho\eta}^B}{y^A}+
G_{\gamma\mu}^A\derpar{G_{\rho\eta}^B}{v_\gamma^A}-
\derpar{G_{\rho\mu}^B}{x^\eta}-F_\eta^A\derpar{G_{\rho\mu}^B}{y^A}-
G_{\gamma\eta}^A\derpar{G_{\rho\mu}^B}{v_\gamma^A}\right)
(\d x^\mu\wedge\d x^\eta )\otimes\derpar{}{v_\rho^B}
\label{curexploc}
\eea
\item
 From the above local expressions
we obtain for a representative multivector field $X$ of the class
$\{ X\}$ associated with the jet field ${\cal Y}$, the expression
\beq
X\equiv\bigwedge_{\mu=1}^m X_\mu=
\bigwedge_{\mu=1}^m \left(\derpar{}{x^\mu}+F_\mu^A\derpar{}{y^A}+
G_{\mu\rho}^A\derpar{}{v_\rho^A}\right)
\label{locmvf}
\eeq
Moreover, as in the case of the bundle $E$,
if $M$ is an orientable manifold and
$\omega\in\df^m(M)$ is the volume form in $M$,
then the above representative can be obtained from the relation
\dst\inn (Y)(\bar\pi^{1^*}\omega )=
\inn\left(\derpar{}{x^1}\wedge\ldots\wedge\derpar{}{x^m}\right)(\bar\pi^{1^*}\omega)\) .
\een

\subsection{Holonomic multivector fields and holonomic jet fields}
\protect\label{hmjfrjb}

We wish to characterize the integrable multivector fields in $J^1E$
whose integral manifolds are canonical prolongations of
sections of $\pi$. In order to achieve this we first recall the same situation for
jet fields \cite{EMR-96}, \cite{Sa-89}.

\begin{definition}
A jet field ${\cal Y}$ is said to be {\sl holonomic}
iff it is integrable and its integral sections $\psi\colon M\to J^1E$ are holonomic; that is,
they are canonical liftings of sections $\phi\colon M\to E$.
\end{definition}

Let ${\bf y}\in J^1J^1E$ with
\dst{\bf y}\stackrel{\pi^1_1}{\mapsto}\bar y\stackrel{\pi^1}{\mapsto}
y\stackrel{\pi}{\mapsto}x\) ,
and $\psi\colon M\to J^1E$ a representative of ${\bf y}$; that is,
${\bf y} =\Tan_x\psi$. Consider the section
$\phi =\pi^1\circ\psi\colon M\to E$, and let
$j^1\phi$ be its canonical prolongation.
Then we are able to define another natural projection
$$
\begin{array}{ccccc}
j^1\pi^1&\colon&J^1J^1E&\longrightarrow&J^1E
\\
& &{\bf y}&\mapsto&j^1(\pi^1\circ\psi )(\bar\pi^1_1({\bf y}))
\\
& &(x^\mu,y^A,v_\mu^A,a_\rho^A,b_{\mu\rho}^A)
&\mapsto&(x^\mu,y^A,a^A_\rho)
\end{array}
$$

\begin{definition}
A jet field ${\cal Y}\colon J^1E\to J^1J^1E$
is said to be a {\rm Second Order Partial Differential Equation (SOPDE)}
or a {\rm semi-holonomic} jet field,
(or also that it verifies the {\rm SOPDE condition}), iff
it is also a section of the projection $j^1\pi^1$; that is,
$j^1\pi^1\circ {\cal Y}={\rm Id}_{J^1E}$.

This is equivalent to saying that ${\cal Y}$ is
a section of the projection $\hat J^2E\to J^1E$,
where $\hat J^2E$ denotes the {\rm semi-holonomic 2-jet manifold}
(see \cite{Sa-89}, p. 173).
\label{sopdejf}
\end{definition}

The relation among integrable, holonomic and SOPDE
jet fields (connections) is the following:

\begin{teor}
A jet field ${\cal Y}\colon J^1E\to J^1J^1E$ is holonomic if, and only if, it is integrable and SOPDE.
\label{holjf}
\end{teor}

The condition $j^1\pi^1\circ {\cal Y}={\rm Id}_{J^1E}$
is locally expressed as follows:
${\cal Y}=(x^\mu ,y^A,v_\mu^A,F_\rho^A,G_{\nu\rho}^A)$
is a SOPDE if, and only if, $F_\rho^A=v_\rho^A$.

For multivector fields in $J^1E$, first we define:

\begin{definition}
A multivector field $X\in\vf^m(J^1E)$ is said to be {\rm holonomic} iff:
\ben
\item
$X$ is integrable.
\item
$X$ is $\bar\pi^1$-transverse.
\item
The integral sections $\psi\colon M\to J^1E$ of $X$ are holonomic.
\een
\end{definition}

In order to study the SOPDE question for multivector fields,
we adopt two different approaches:
The first one, which is similar to the characterization
of {\sl second order} (that is, holonomic) vector fields
in time-dependent mechanics, is based on the study of the two
different bundle structures of $\Tan J^1E$ over $\Tan E$.
The second one uses the {\sl structure canonical form}
$\theta\in\df^1(J^1E,\pi^{1^*}{\rm V}(\pi))$ (see the appendix).

In the first approach, we have a natural vector bundle projection
$\Tan\pi^1\colon\Tan J^1E\to\Tan E$
and another one $\kappa\colon\Tan J^1E\to\Tan E$ defined by
$$
\kappa (\bar y,\bar u):=
\Tan_{\bar\pi^1(\bar y)}\phi (\Tan_{\bar y}\bar\pi^1(\bar u))
$$
where $(\bar y,\bar u)\in\Tan J^1E$ and $\phi\in\bar y$.
If $(W;x^\mu,y^A,v^A_\mu)$ is a local natural chart in $J^1E$
and $\bar y\in J^1E$ with
$\bar y\stackrel{\pi^1}{\to}y\stackrel{\pi}{\to}x$, then
$$
\kappa\left( \alpha^\mu\derpar{}{x^\mu}\Big\vert_{\bar y}+
\beta^A\derpar{}{y^A}\Big\vert_{\bar y}+
\lambda^A_\mu\derpar{}{v^A_\mu}\Big\vert_{\bar y}\right)=
\alpha^\mu\derpar{}{x^\mu}\Big\vert_y+
v^A_\mu(\bar y)\alpha^\mu\derpar{}{y^A}\Big\vert_y
$$

This projection is extended in a natural way to $\Lambda^m\Tan J^1E$,
and so we have the following diagram
$$
\begin{array}{cccccccccc}
& & & \tau_{J^1E} & & & \Lambda^m\tau_{J^1E} & & &
\\
& \Tan J^1E& &
\begin{picture}(60,10)(0,0)
\put(0,5){\vector(1,0){50}}
\end{picture}
& J^1E & &
\begin{picture}(50,10)(0,0)
\put(50,7){\vector(-1,0){50}}
\put(0,3){\vector(1,0){50}}
\end{picture}
& &\Lambda^m\Tan J^1E&
\\
& & & & & & X & & &
\\
\begin{picture}(5,50)
\put(0,23){\mbox{$\kappa$}}
\end{picture}
&
\begin{picture}(10,50)(0,0)
\put(3,50){\vector(0,-1){50}}
\end{picture}
\begin{picture}(10,50)(0,0)
\put(7,50){\vector(0,-1){50}}
\end{picture}
&
\begin{picture}(5,50)
\put(0,23){\mbox{$\Tan\pi^1$}}
\end{picture}
& &
\begin{picture}(2,50)(0,0)
\put(1,50){\vector(0,-1){50}}
\end{picture}
&
\begin{picture}(5,50)
\put(0,23){\mbox{$\pi^1$}}
\end{picture}
& &
\begin{picture}(15,50)
\put(0,23){\mbox{$\Lambda^m\Tan\pi^1$}}
\end{picture}
&
\begin{picture}(10,50)(0,0)
\put(3,50){\vector(0,-1){50}}
\end{picture}
\begin{picture}(10,50)(0,0)
\put(7,50){\vector(0,-1){50}}
\end{picture}
&
\begin{picture}(5,50)
\put(-15,23){\mbox{$\Lambda^m\kappa$}}
\end{picture}
\\
& & & & & & & & &
\\
& \Tan E& &
\begin{picture}(50,10)(0,0)
\put(0,5){\vector(1,0){50}}
\end{picture}
& E & &
\begin{picture}(50,10)(0,0)
\put(50,5){\vector(-1,0){50}}
\end{picture}
& &\Lambda^m\Tan E&
\\
& & & \tau_E & & & \Lambda^m\tau_E & & &
\end{array}
$$
and we can define:

\begin{definition}
A $\bar\pi^1$-transverse multivector field $X\in\vf^m(J^1E)$
verifies the {\rm SOPDE condition}
(we also say that it is a {\rm SOPDE} or
a {\rm semi-holonomic} multivector field) iff
$\Lambda^m\kappa\circ X=\Lambda^m\Tan\pi^1\circ X$.
\label{sopdemvdef}
\end{definition}

The condition for $X$ to be $\bar\pi^1$-transverse
is necessary in order to relate this definition to the fact that,
if $X$ is integrable, then its integral sections are canonical prolongations
of sections of $\pi$ (as will be seen later).

With regard to the second approach, observe that for every $X\in\vf (J^1E)$,
we obtain that $\theta (X)\equiv\inn(\theta)X$ is
a section of the bundle $\pi^{1^*}{\rm V}(\pi)\to J^1E$
(that is, an element of $\Gamma(\pi^{1^*}{\rm V}(\pi))$).
Then we have a natural extension of the action of $\theta$ to
every $X\in\vf^m(J^1E)$.
Furthermore, denoting ${\rm V}^m(\pi):={\rm ker}\,\Lambda^m\Tan\pi$
and recalling the definition of $\theta$
(see definition  \ref{theta} in the appendix)
we can define:

\begin{definition}
The {\rm structure canonical $m$-form} of $J^1E$ is the $m$-form $\theta^m$
in $J^1E$ with values on $\pi^{1^*}{\rm V}^m(\pi)$, which is defined by
$$
(\theta^m (X))(\bar y)\equiv (\inn(\theta^m)X)(\bar y):=
(\Lambda^m\Tan_{\bar y}\pi^1-
\Lambda^m\Tan_{\bar y}(\phi\circ\bar\pi^1))(X_{\bar y})
$$
where $X\in\vf^m(J^1E)$, $\bar y\in J^1E$
and $\phi$ is a representative of $\bar y$.
\label{thetam}
\end{definition}

Observe that $\theta^m\in\df^m(J^1E,\pi^{1^*}{\rm V}^m(\pi))$
and that \dst\theta^m\not=\overbrace{\theta\wedge\ldots\wedge\theta}^{m}\) .

Now we can state:

\begin{prop}
Let $X\in\vf^m(J^1E)$ be $\bar\pi^1$-transverse
and locally decomposable. Then the following conditions are equivalent:
\ben
\item
$X$ is a SOPDE.
\item
 $\inn(\theta^m)X=0$.
\item
 $\inn(\theta )X=0$.
\item
If $(W;x^\mu,y^A,v^A_\mu)$ is a natural chart in $J^1E$,
then the local expression of $X$ is
\beq
X\equiv\bigwedge_{\mu=1}^m X_\mu=
\bigwedge_{\mu=1}^m f_\mu\left(\derpar{}{x^\mu}+v_\mu^A\derpar{}{y^A}+
G_{\mu\rho}^A\derpar{}{v_\rho^A}\right)
\label{locsopdegen}
\eeq
where $f_\mu$ are non-vanishing functions.
\item
$\inn (\alpha)X=0$, for every $\alpha\in{\cal M}_c$;
where ${\cal M}_c$ denotes the {\sl contact module} in $J^1E$
(see the appendix).
\een
\label{canform}
\end{prop}
\proof
Let $(W;x^\mu,y^A,v^A_\mu)$ be a natural chart in $J^1E$.
Suppose that $X\in\vf^m(J^1E)$ is $\bar\pi^1$-transverse
and locally decomposable, then in this chart we have
$$
X=\bigwedge_{\mu=1}^m f_\mu\left(\derpar{}{x^\mu}+F_\mu^A\derpar{}{y^A}+
G_{\mu\rho}^A\derpar{}{v_\rho^A}\right)
$$
with $f_\mu$ non-vanishing functions.

\quad ( 1 $\Leftrightarrow$ 2 ) \quad
It is obvious.

\quad ( 2 $\Leftrightarrow$ 4 ) \quad
If $\bar y\in U\subset J^1E$ and $\phi\in\bar y$, then
\beann
(\theta^m (X))(\bar y) &=&
(\Lambda^m\Tan_{\bar y}\pi^1
-\Lambda^m\Tan_{\bar y}(\phi\circ\bar\pi^1))(X_{\bar y})
\\ &=&
\bigwedge_{\mu=1}^m f_\mu (\bar y)\left(\derpar{}{x^\mu}+
F_\mu^A (\bar y)\derpar{}{y^A}\right)_y-
\bigwedge_{\mu=1}^m f_\mu (\bar y)\left(\derpar{}{x^\mu}+
v_\mu^A (\bar y)\derpar{}{y^A}\right)_y
\eeann
therefore $\theta^m (X)=0$ $\Leftrightarrow$ $F^A_\mu=v^A_\mu$,
for every $A,\mu$.

\quad ( 3 $\Leftrightarrow$ 4 ) \quad
Since
\dst\theta =\left(\d y^B-v^B_\rho\d x^\rho\right)\otimes\derpar{}{y^B}\) ,
we have that
\beann
\inn(\theta)X &=& \left[
f^1(F^B_1-v^B_1)\bigwedge_{\mu\not= 1}f_\mu
\left(\derpar{}{x^\mu}+
F_\mu^A\derpar{}{y^A}+G_{\mu\rho}^A\derpar{}{v_\rho^A}\right) - \right.
\\ & &
\left. \ f^2(F^B_2-v^B_2)\bigwedge_{\mu\not= 2}
f_\mu\left(\derpar{}{x^\mu}+F_\mu^A\derpar{}{y^A}+
G_{\mu\rho}^A\derpar{}{v_\rho^A}\right) + \ldots \right]
\otimes\derpar{}{y^B}
\eeann
therefore $\inn (X)\theta=0$ $\Leftrightarrow$ $F^A_\mu=v^A_\mu$,
for every $A,\mu$.

\quad ( 4 $\Leftrightarrow$ 5 ) \quad
It is immediate from the item 4, and bearing in mind
the local expressions given in the appendix.
\qed

The relation between integrable, holonomic and SOPDE multivector fields in $J^1E$ is:

\begin{teor}
A multivector field $X\in\vf^m(J^1E)$ is holonomic if, and only if, it is integrable and SOPDE.
\label{holmvf}
\end{teor}
\proof
Let $\bar y\in J^1E$ and $\psi\colon M\to J^1E$ an integral section
of $X$ at $\bar y$, with $\psi (x)=\bar y$.
For $Z\in\Lambda^m\Tan_xM$ with $Z_x\not= 0$,
as $\psi$ is an integral section, there exists
$\lambda\in\Real$ with $\lambda\not= 0$
such that $\Lambda^m\Tan_x\psi(Z_x)=\lambda X_{\bar y}$
and hence
$$
(\inn(\theta)X)(\bar y)=
\frac{1}{\lambda}\inn(\theta)(\Lambda^m\Tan_x\psi(Z_x))=
\frac{1}{\lambda}\inn(\psi^*\theta)Z_x=0
$$
therefore $\inn (X)\theta=0$ $\Leftrightarrow$ $\psi^*\theta=0$,
which is the necessary and sufficient condition for $\psi$ to be holonomic
(see proposition \ref{holsec} in the appendix).
\qed

Finally, as a consequence of
theorems \ref{vfmvf1}, \ref{holjf} and \ref{holmvf} we have:

\begin{teor}
There is a bijective correspondence between the set of
orientable jet fields ${\cal Y}\colon J^1E\to J^1J^1E$
(i. e., the set of orientable connections $\nabla$ in $\bar\pi^1\colon J^1E\to M$)
and the set of the equivalence classes
of locally decomposable and $\bar\pi^1$-transverse multivector fields
$\{ X\}\subset\vf^m(J^1E)$.
They are characterized by the fact that
${\cal D}({\cal Y})={\cal D}(X)$.
In addition:
\ben
\item
The jet field ${\cal Y}$ is integrable if, and only if,
so is $X$, for every $X\in\{ X\}$.
\item
The jet field ${\cal Y}$ is a SOPDE if, and only if,
so is $X$, for every $X\in\{ X\}$.
\item
The jet field ${\cal Y}$ is holonomic if, and only if,
so is $X$, for every $X\in\{ X\}$.
\een
\label{vfmvf}
\end{teor}

If ${\cal Y}$ is a SOPDE, then we can choose (see the local expressions (\ref{locmvf}) and (\ref{locsopdegen}))
\beq
X\equiv
\bigwedge_{\mu=1}^m \left(\derpar{}{x^\mu}+v_\mu^A\derpar{}{y^A}+
G_{\mu\rho}^A\derpar{}{v_\rho^A}\right)
\label{locsopde}
\eeq
as a representative multivector field of $\{ X\}$.
For a SOPDE multivector field or jet field, on the other hand, if
a section \dst j^1\phi =\left( x^\mu ,f^A,\derpar{f^A}{x^\mu}\right)\)
has to be an integral section of such a field, the necessary and sufficient condition is that
$\phi$ is the solution of the following system of (second order) partial differential equations
\beq
G_{\nu\rho}^A\left(x^\mu ,f^A,\derpar{f^A}{x^\mu}\right) =
\frac{\partial^2f^A}{\partial x^\rho\partial x^\nu}
\label{sisteq2}
\eeq
which justifies the nomenclature.  It is important to remark that,
since the integrability of a class of multivector fields is equivalent to demanding that
the curvature ${\cal R}$ of the connection associated with this class
vanishes everywhere; the system (\ref{sisteq2}) has solution
if, and only if, the following additional system of equations
(which arise from the condition (\ref{curexploc})  for a SOPDE) holds
(for every $B,\mu,\rho,\eta$)
\beq
\begin{array}{lll}
0 &=& G^B_{\eta\mu}-G^B_{\mu\eta}
\\
0 &=& \derpar{G_{\eta\rho}^B}{x^\mu}+v_\mu^A\derpar{G_{\eta\rho}^B}{y^A}+
G_{\mu\gamma}^A\derpar{G_{\eta\rho}^B}{v_\gamma^A}-
\derpar{G_{\mu\rho}^B}{x^\eta}-v_\eta^A\derpar{G_{\mu\rho}^B}{y^A}-
G_{\eta\gamma}^A\derpar{G_{\mu\rho}^B}{v_\gamma^A}
\end{array}
\label{curvcero}
\eeq
Observe that this is a system of \dst\frac{1}{2}Nm(m-1)\) linear relations
(the first group of equations) and $\frac{1}{2}Nm^2(m-1)\)
partial differential equations (the second group).

\section{Application to Lagrangian classical field theories}

At this point, our goal is to show that
the Lagrangian formalism for field theories can also be established
using jet fields in $J^1J^1E$, their associated connections
in $J^1E$ or, equivalently, multivector fields in $J^1E$.
Then, we use multivector fields for discussing the
main features of the evolution equations.

\subsection{Background on Lagrangian field theories}
\protect\label{eblft}

Henceforth we assume that $M$ is a $m$-dimensional oriented manifold and
$\omega\in\df^m(M)$ is the volume $m$-form on $M$.

A {\sl classical field theory} is described by its {\sl configuration bundle}
$\pi\colon E\to M$ (where $M$ is an oriented manifold with volume form $\omega$),
and a {\sl Lagrangian density} which is
a $\bar\pi^1$-semibasic $m$-form on $J^1E$.
A Lagrangian density is usually written as
$\Lag =\lag (\bar\pi^{1^*}\omega)$, where $\lag\in\Cinfty (J^1E)$
is the {\sl Lagrangian function} associated with $\Lag$ and $\omega$.
In a natural system of coordinates this expression is
$\Lag = \lag (x^\mu ,y^A,v^A_\mu )\d x^1\wedge\ldots\wedge\d x^m$.
Then a {\sl Lagrangian system} is a couple $\ls$.
The states of the field are the sections of $\pi$ which are critical for the functional
$$
\begin{array}{ccccc}
{\bf L}&\colon&\Gamma_c(M,E)&\longrightarrow&\Real
\\
& &\phi&\mapsto&\int_M(j^1\phi)^*\Lag
\end{array}
$$
where $\Gamma_c(M,E)$ is the set of compact supported sections of $\pi$.

In order to characterize these critical sections, and
using the {\sl vertical endomorphism} ${\cal V}$ of the bundle $J^1E$
(see the appendix), we construct the differentiable forms
$$
\Theta_{\Lag}:=\inn({\cal V})\Lag+\Lag\in\df^{m}(J^1E)
\quad ;\quad
\Omega_{\Lag}:= -\d\Theta_{\Lag}\in\df^{m+1}(J^1E)
$$
which are called the {\sl Poincar\'e-Cartan $m$ and $(m+1)$-forms}
associated with the Lagrangian density $\Lag$.
The Lagrangian system is said to be {\sl regular} iff
$\Omega_{\Lag}$ is $1$-nondegenerate and, as a consequence,
$(J^1E,\Omega_{\Lag})$ is a multisymplectic manifold \cite{CIL-96a}.

In a natural chart $(x^\mu ,y^A,v^A_\mu )$ in $J^1E$,
the expressions of the above forms are:
\bea
\Theta_{\Lag}&=&
\derpar{\lag}{v^A_\mu}\d y^A\wedge\d^{m-1}x_\mu -
\left(\derpar{\lag}{v^A_\mu}v^A_\mu -\lag\right)\d^mx
\nonumber \\
\Omega_{\Lag}&=&
-\frac{\partial^2\lag}{\partial v^B_\nu\partial v^A_\mu}
\d v^B_\nu\wedge\d y^A\wedge\d^{m-1}x_\mu 
\nonumber \\ & &
-\frac{\partial^2\lag}{\partial y^B\partial v^A_\mu}\d y^B\wedge
\d y^A\wedge\d^{m-1}x_\mu +
\frac{\partial^2\lag}{\partial v^B_\nu\partial v^A_\mu}v^A_\mu
\d v^B_\nu\wedge\d^mx  +
\nonumber \\ & &
\left(\frac{\partial^2\lag}{\partial y^B\partial v^A_\mu}v^A_\mu -\derpar{\lag}{y^B}+
\frac{\partial^2\lag}{\partial x^\mu\partial v^B_\mu}
\right)\d y^B\wedge\d^mx
\label{omegalag}
\eea
(where 
\dst\d^{m-1}x_\mu\equiv\inn\left(\derpar{}{x^\mu}\right)\d^mx\) )
and the condition of regularity is equivalent to demand that
\dst det\left(\frac{\partial^2\lag}{\partial v^A_\mu\partial v^B_\nu}(\bar y)\right)\not= 0\) ,
for every $\bar y\in J^1E$.

Then, the (compact-supported) critical sections $\phi\colon M\to E$
of the variational problem can be characterized as follows:

\begin{prop}
Let $\ls$ be a Lagrangian system.
The critical sections of the variational problem posed by $\Lag$
are those which satisfy the following equivalent conditions
\ben
\item
\dst\frac{\d}{\d t}\Big\vert_{t=0}\int_M(j^1\phi_t)^*\Lag = 0\) ,
being $\phi_t =\tau_t\circ \phi$, where $\{\tau_t\}$
is a local one-parameter group of any
$\pi$-vertical and compact supported vector field $Z\in\vf (E)$.
\item
\dst(j^1\phi)^*\inn (X)\Omega_{\Lag} = 0\) ,
for every $X\in\vf (J^1E)$.
\item
The local components of $\phi$ satisfy the {\sl Euler-Lagrange equations}:
$$
\left(\derpar{\lag}{y^A}\right)_{j^1\phi}-
\derpar{}{x^\mu}\left(\derpar{\lag}{v_\mu^A}\right)_{j^1\phi} = 0
$$
\een
\label{carcrisec}
\end{prop}

(See, for instance, \cite{EMR-96} for the proof,
and also \cite{BSF-88}, \cite{Gc-73}, \cite{Sd-94c}, \cite{Sa-87}
as complementary references on all these topics):

\subsection{The evolution equations in terms of connections and multivector fields}
\protect\label{lfcft}

As is well known, the evolution equations in analytical mechanics can be written
in a geometric-algebraic manner using vector fields,
and then the solutions of the variational problem can be obtained
as integral curves of these vector fields.
Now we wish to do the same with the evolution equations of field theories,
using jet fields (connections) and their equivalent multivector fields,
and then obtaining the critical sections as integral sections of these elements.

As a previous step, we must define the action of jet fields on
differential forms (see \cite{EMR-96}).

Let ${\cal Y}\colon J^1E\to J^1J^1E$ be a jet field.
A map $\bar{\cal Y}\colon\vf (M)\to\vf (J^1E)$
can be defined in the following way:
let $Z\in\vf (M)$, then
$\bar{\cal Y}(Z)\in\vf(J^1E)$ is the vector field defined as
$$
\bar {\cal Y}(Z)(\bar y):=
(\Tan_{\bar\pi^1 (\bar y)}\psi )(Z_{\bar\pi^1 (\bar y)})
$$
for every $\bar y\in J^1E$ and $\psi\in {\cal Y}(\bar y)$.
This map is an element of $\df^1(M)\otimes_M\vf (J^1E)$
and its local expression is
$$
\bar{\cal Y}\left( f^\mu\derpar{}{x^\mu}\right)=
f^\mu\left(\derpar{}{x^\mu}+F_\mu^A\derpar{}{y^A}+
G_{\rho\mu}^A\derpar{}{v_\rho^A}\right)
$$
This map induces an action of ${\cal Y}$ on
$\df (J^1E)$. In fact, let $\xi\in\df^{m+j}(J^1E)$, and
with $j\geq 0$, we define
$\inn ({\cal Y})\xi\colon\vf (M)\times\stackrel{(m)}{\ldots}\times\vf (M)
\longrightarrow\df^j(J^1E)$
given by
$$
((\inn ({\cal Y})\xi )(Z_1,\ldots ,Z_m))(\bar y;X_1,\ldots ,X_j):=
\xi (\bar y;\bar{\cal Y}(Z_1),\ldots ,\bar{\cal Y}(Z_m),X_1,\ldots ,X_j)
$$
for $Z_1,\ldots ,Z_m\in\vf (M)$ and $X_1,\ldots ,X_j\in\vf (J^1E)$.
This is a $\Cinfty (M)$-linear and alternate map on the vector fields $Z_1,\ldots ,Z_m$.
This map $\inn ({\cal Y})\xi$ so-defined, extended by zero to forms of degree $p<m$,
is called the {\sl inner contraction} of the jet field
${\cal Y}$ and the differential form $\xi$.

The following result characterizes the action
of integrable jet fields on forms
(see \cite{EMR-96} for the proof):

\begin{prop}
Let ${\cal Y}$ be an integrable jet field and
$\xi\in\df^p(J^1E)$ with $p\geq m$.
The necessary and sufficient condition for the integral sections of ${\cal Y}$,
$\psi\colon M\to J^1E$, to verify the relation
$\psi^*\inn (X)\xi =0$, for every $X\in\vf (J^1E)$,
is $\inn ({\cal Y})\xi =0$.
\label{tecnic}
\end{prop}

At this point the following theorem gives the characterization of the critical sections
of a Lagrangian field theory in terms of connections (jet fields) and multivector fields.

\begin{teor}
Let $\ls$ be a Lagrangian system.
The critical sections of the variational problem posed by $\Lag$
satisfy the following equivalent conditions:
\ben
\item
They are the integral sections of a holonomic jet field
${\cal Y}_{\Lag}\colon J^1E\to J^1J^1E$, such that
$$
\inn ({\cal Y}_{\Lag})\Omega_{\Lag}=0
$$
\item
They are the integral sections of a holonomic jet field
${\cal Y}_{\Lag}\colon J^1E\to J^1J^1E$, such that
$$
\inn (\nabla_{\Lag})\Omega_{\Lag}=(m-1)\Omega_{\Lag}
$$
where $\nabla_{\Lag}$ is the associated connection form.
\item
They are the integral sections of a class of
holonomic multivector fields  $\{ X_{\Lag}\}\subset\vf^m(J^1E)$, such that
$$
\inn (X_{\Lag})\Omega_{\Lag}=0
$$
\een
\label{important}
\end{teor}
\proof
These results are consequences of all the above results.
In particular:
\ben
\item
It follows from the second item in proposition \ref{carcrisec} and proposition \ref{tecnic}.
\item
See \cite{LMM-95} and \cite{Sa-87} for this proof.
\item
It is a consequence of the first item and theorem \ref{vfmvf}.
\een
\qed

The conditions in this theorem are the version of the Euler-Lagrange equations in terms of
jet fields, connection forms and multivector fields respectively.

\subsection{Analysis of the evolution equations for regular Lagrangians}
\protect\label{aeerlft}

Using these formulations and, in particular, the multivector fields
formalism, it is possible to explore some important properties
of the evolution equations in field theories,
specifically the existence and multiplicity of solutions.
As a first situation, we can consider the case of regular Lagrangian systems.

First, let us recall that the equivalent problem in
(time-dependent) mechanics consists in finding vector fields
$X_{\Lag}\in\vf(\Tan Q\times\Real)$
(where $\Tan Q$ is the velocity phase space of the system),
such that they verify the dynamical equation $\inn(X_{\Lag})\Omega_{\Lag}=0$,
and they are holonomic vector fields (their integral curves
are canonical lifting of curves in $Q$; that is, sections
of the projection $Q\times\Real\to\Real$).
In this case, the regularity of the Lagrangian function
assures the existence and uniqueness of an holonomic vector field
which is the solution of the intrinsic equation of motion,
with the additional condition $\inn(X_{\Lag})\d t=1$.

In field theories we search for (classes of) non-vanishing and locally decomposable multivector fields $X_{\Lag}\in\vf^m(J^1E)$ such that:
\beq
\begin{array}{lll}
1. & \mbox{The equation $\inn (X_{\Lag})\Omega_{\Lag}=0$ holds.}
& \qquad\qquad\qquad\qquad\qquad\qquad\qquad\qquad\qquad\qquad\qquad
\\
2. & \mbox{$X_{\Lag}$ are SOPDE.}
& \qquad\qquad\qquad\qquad\qquad\qquad\qquad\qquad\qquad\qquad\qquad
\\
3. & \mbox{$X_{\Lag}$ are integrable.}
& \qquad\qquad\qquad\qquad\qquad\qquad\qquad\qquad\qquad\qquad\qquad
\end{array}
\label{conditions}
\eeq

{\bf Remarks}:
\bit
\item
If a multivector field satisfies condition 2 and also verifies that
$\inn (X_{\Lag})(\bar\pi^{1*}\omega)=1$, then its local expression is (\ref{locmvf}).
\item
Let us recall that, if $\{ X_{\Lag}\}\subset\vf^m(J^1E)$ is a class of
locally decomposable and $\bar\pi^1$-transverse multivector fields,
then condition 3 is equivalent to demanding
that the curvature ${\cal R}$ of the connection associated with this class
vanishes everywhere.
\eit

Following the terminology introduced in \cite{LMM-95} and \cite{Sa-89}
for connections, we define:

\begin{definition}
$X_{\Lag}\in\vf^m(J^1E)$ is said to be an
{\rm Euler-Lagrange multivector field} for $\Lag$ iff
it is locally decomposable and verifies conditions 1 and 2 of \ref{conditions}.
\end{definition}

Regularity of Lagrangians assures the existence of Euler-Lagrange
multivector fields (or connections). In fact:

\begin{teor}
{\rm (Existence and local multiplicity of Euler-Lagrange solutions)}
Let $\ls$ be a regular Lagrangian system.
\ben
\item
There exist classes of locally decomposable
and $\bar\pi^1$-transverse multivector fields
$\{ X_{\Lag}\}\subset\vf^m(J^1E)$
(and hence equivalent jet fields ${\cal Y}_{\Lag}\colon J^1E\to J^1J^1E$
with associated connection forms $\nabla_{\Lag}$), such that
\ben
\item
The following equivalent conditions hold
\beq
\inn (X_{\Lag})\Omega_{\Lag}=0
\quad ,\quad
\inn ({\cal Y}_{\Lag})\Omega_{\Lag}=0
\quad ,\quad
\inn (\nabla_{\Lag})\Omega_{\Lag}=(m-1)\Omega_{\Lag}
\label{evoleq}
\eeq
\item
$X_{\Lag}\in\{ X_{\Lag}\}$ are SOPDE (and therefore so is the corresponding
${\cal Y}_{\Lag}$).
\een
\item
In a local system the above solutions depend on $N(m^2-1)$ arbitrary functions.
\een
\label{holsecreg}
\end{teor}
\proof
In order to prove this, we use the multivector field formalism
(for another proof of these statement using the Ehresmann connection formalism see \cite{LMM-95} and \cite{Sa-89}).
\ben
\item
First we analyze the local existence of solutions and then their global extension.

In a chart of natural coordinates in $J^1E$,
the expression of $\Omega_{\Lag}$ is (\ref{omegalag}),
and taking the multivector field given in (\ref{locmvf}) as the
representative of the class $\{ X_{\Lag}\}$, from the relation
$\inn (X_{\Lag})\Omega_{\Lag}=0$ we obtain the following conditions:
\bit
\item
The coefficients on $\d v^A_\mu$ must vanish:
\beq
0 = (F^B_\mu-v^B_\mu)\frac{\partial^2\lag}{\partial v^A_\nu\partial v^B_\mu}
\qquad (\mbox{for every $A,\nu$})
\label{eqsG1}
\eeq
But, if $\Lag$ is regular, the matrix
\dst\left(\frac{\partial^2\lag}{\partial v^A_\nu\partial v^B_\mu}\right)\)
is regular. Therefore $F^B_\mu=v^B_\mu$ (for every $B,\mu$);
which proves that if $X_{\Lag}$ exists it is a SOPDE.
\item
The coefficients on $\d y^A$ must vanish
\beq
0 = \derpar{\lag}{y^A}-\frac{\partial^2\lag}{\partial x^\mu\partial  v^A_\mu}-
\frac{\partial^2\lag}{\partial y^B\partial v^A_\mu}F^B_\mu-
\frac{\partial^2\lag}{\partial v^B_\nu\partial  v^A_\mu}G^B_{\nu\mu}+
\frac{\partial^2\lag}{\partial y^A\partial  v^B_\mu}(F^B_\mu-v^B_\mu)
\quad (A=1,\ldots ,N)
\label{eqsG2}
\eeq
and taking into account that we have obtained $F^B_\mu=v^B_\mu$, these relations lead to
\beq
\frac{\partial^2\lag}{\partial v^B_\nu\partial  v^A_\mu}G^B_{\nu\mu}=
\derpar{\lag}{y^A}-\frac{\partial^2\lag}{\partial x^\mu\partial  v^A_\mu}-
\frac{\partial^2\lag}{\partial y^B\partial v^A_\mu}v^B_\mu
\qquad (A=1,\ldots ,N)
\label{eqsG}
\eeq
which is a system of $N$ linear equations on the functions $G^B_{\nu\mu}$.
This is a compatible system as a consequence of the regularity of $\Lag$,
since the matrix of the coefficients has (constant) rank equal to $N$
(observe that the matrix of this system is obtained as a rearrangement
of rows of the Hessian matrix).
\item
Finally, from the above conditions, we obtain that the coefficients on
$\d x^\mu$ vanish identically.
\eit
These results allow us to assure the local existence of (classes of) multivector fields
satisfying the desired conditions. The corresponding global solutions
are then obtained using a partition of unity subordinated
to a covering of $J^1E$ made of local natural charts.

(Observe that, if
\dst j^1\phi =\left( x^\mu ,f^A,\derpar{f^A}{x^\mu}\right)\)
is an integral section of $X_{\Lag}$, then
\dst v^A_\mu=\derpar{f^A}{x^\mu}\) and
\dst G_{\nu\mu}^A=\frac{\partial^2f^A}{\partial x^\nu\partial x^\mu}\) ,
and therefore the equations (\ref{eqsG}) are the Euler-Lagrange equations
for the section $\phi$).
\item
In a chart of natural coordinates in $J^1E$,
the expression of a SOPDE multivector field
$X_{\Lag}\in\{ X_{\Lag}\}$ is given by (\ref{locsopde}).
So, it is determined by the $Nm^2$ coefficients $G^B_{\nu\mu}$,
which are related  by the $N$ independent equations (\ref{eqsG}).
Therefore, there are $N(m^2-1)$ arbitrary functions.
\een
\qed

{\bf Comments}:
\bit
\item
If, in addition, $\{ X_{\Lag}\}$ is integrable;
then every $X_{\Lag}\in\{ X_{\Lag}\}$ is a holonomic multivector field
(and, therefore, ${\cal Y}_{\Lag}$ is also a holonomic jet field).
\item
As additional conditions on the functions $G^A_{\rho\mu}$
must be imposed in order to assure that
$X_{\Lag}$ is integrable, hence the number of arbitrary functions
will be in general less than $N(m^2-1)$ (see the next section).
\eit

The last problem consists in studying if it is possible to find
a class of integrable Euler-Lagrange multivector fields.
This question is posed in the following terms:
taking the suitable representatives in each class, the set of
Euler-Lagrange multivector fields $X_{\Lag}$ can be locally written in the form
(\ref{locsopde}), where the coefficients $G^A_{\mu\nu}$
verify the equations (\ref{eqsG}). This is a linear system and,
then, if ${\cal G}^A_{\mu\nu}$ is a particular solution of this system,
we can write
$$
X_{\Lag}=
\bigwedge_{\mu=1}^m \left(\derpar{}{x^\mu}+v_\mu^A\derpar{}{y^A}+
{\cal G}_{\mu\gamma}^B\derpar{}{v_\gamma^B}+g_{\mu\eta}^C\derpar{}{v_\eta^C}\right)
$$
where $g^C_{\eta\mu}$ satisfy the homogeneous linear system
$$
\frac{\partial^2\lag}{\partial v^C_\eta\partial v^A_\mu}g^C_{\eta\mu}=0
\qquad (A=1,\ldots ,N)
$$
As the problem now is to select the integrable Euler-Lagrange multivector fields, from
the solutions of this system we can choose those that make
the corresponding multivector field $X_{\Lag}$ verify the integrability condition;
that is, the curvature of the associated connection $\nabla_{\Lag}$
vanishes (equations (\ref{curvcero})).

As far as we know,
since we have a system of partial differential equations with linear restrictions,
there is no way for assuring the existence of an integrable solution;
or how to select it;
Observe that, considering the Euler-Lagrange equations for the coefficients $G^A_{\mu\nu}$ (equations (\ref{eqsG})), together with the integrability conditions
(equations (\ref{curvcero})), we have
\dst N+\frac{1}{2}Nm(m-1)\) linear equations and 
\dst\frac{1}{2}Nm^2(m-1)\) partial differential equations.
Then, if the set of linear restrictions (made of the equations (\ref{eqsG})
and the first group of equations (\ref{curvcero})) allow us to isolate
\dst N+\frac{1}{2}Nm(m-1)\) coefficients $G^A_{\mu\nu}$ as functions on the remaining ones;
and the set of \dst\frac{1}{2}Nm^2(m-1)\) partial differential equations
(the second group of equations (\ref{curvcero})) on these remaining coefficients
satisfies the conditions on {\sl Cauchy-Kowalewska's theorem} \cite{Di-74},
then the existence of integrable Euler-Lagrange multivector fields is assured
(see the example in section \ref{example}).

On the other hand, if we can eventually select some particular
Euler-Lagrange multivector field solution,
then the integrability algorithm developed in section \ref{ia} can be applied in order to find a submanifold where this multivector field is integrable (if it exists).

\subsection{Discussion on the evolution equations for singular Lagrangian field theories}
\protect\label{dslft}

Next we discuss the case of
singular Lagrangian systems ($\Omega_{\Lag}$ is 1-degenerate).
As in the regular case, the existence of
multivector fields (or connections) verifying conditions (\ref{conditions}) is not assured.
The difference with the regular case lies in the following facts:
\ben
\item
A locally decomposable and $\bar\pi^1$-transverse multivector field
$X\in\vf^m(J^1E)$ verifying condition 1on $J^1E$
(the evolution equation) does not necessarily exist.
\item
If such a multivector field exists, it does not necessarily hold condition 2
(SOPDE condition).
\een
As in the regular case, the existence of Euler-Lagrange multivector fields
does not imply their integrability.
Nevertheless, it is possible for these integrable multivector fields to exist
on a submanifold of $J^1E$. So we can state the following problem:
to look for a submanifold $S\hookrightarrow J^1E$ where
integrable Euler-Lagrange multivector fields $X_\Lag\in\vf^m(J^1E)$ exist;
and then their integral sections are contained in $S$.
Observe that these conditions imply that $\bar\pi^1\vert_S\colon S\to M$
is onto on $M$, because $X_{\Lag}$ is SOPDE.

As a first step, we do not consider the integrability condition.
(In this way, the situation is quite analogous to that in (time-dependent) mechanics
\cite{CLM-94}, \cite{LMM-96}).

The procedure is algorithmic
(from now on we suppose that all the multivector fields are locally decomposable):
\bit
\item
First, let $S_1$ be the set of points of $J^1E$
where Euler-Lagrange multivector fields do exist; that is
$$
S_1:=\{ \bar y\in J^1E \ ;\ \exists X_{\Lag}\in\vf^m(J^1E)\
\mbox{such that}
\left\{
\begin{array}{c}
(\inn(X_{\Lag})\Omega_{\Lag})(\bar y)=0
\\
(\inn(X_{\Lag})(\bar\pi^{1*}\omega))(\bar y)=1
\\
\mbox{$X_{\Lag}$ is a SOPDE at $\bar y$}
\end{array}
\right\}
\}
$$
We assume that $S_1$ is a non-empty (closed) submanifold of $J^1E$.

This is the {\sl compatibility condition}.
\item
Now, denote by $\vf^m_{\Lag}(J^1E,S_1)$ the set of multivector fields
in $J^1E$ which are Euler-Lagrange on $S_1$.
Let $X_{\Lag}\colon S_1\to\Lambda^m\Tan J^1E\vert_{S_1}$ be in $\vf^m_{\Lag}(J^1E,S_1)$.
If, in addition, $X_{\Lag}\colon S_1\to\Lambda^m\Tan S_1$; that is, $X_{\Lag}\in\vf^m(S_1)$,
then we say that $X_{\Lag}$ is a solution on $S_1$.
Nevertheless, this last condition is not assured except perhaps in a set of points
$S_2\subset S_1\subset J^1E$,
which we will assume to be a (closed) submanifold, and which is defined by
$$
S_2:=\{ \bar y\in S_1 \ ;\ \exists X_{\Lag}\in\vf^m_{\Lag}(J^1E,S_1)\
\mbox{such that}\
X_{\Lag}(\bar y)\in\Lambda^m\Tan_{\bar y}S_1\}
$$
This is the so-called {\sl consistency} or {\sl tangency condition}.
\item
In this way, a sequence of (closed) submanifolds,
$\ldots \subset S_i\subset\ldots \subset S_1\subset J^1E$,
is assumed to be obtained, each one of them being defined as
$$
S_i:=\{ \bar y\in S_{i-1} \ ;\ \exists X_{\Lag}\in\vf^m_{\Lag}(J^1E,S_{i-1})\
\mbox{such that}\
X_{\Lag}(\bar y)\in\Lambda^m\Tan_{\bar y}S_{i-1}\}
$$
It is important to point out that, for every submanifold $S_i$ of this sequence,
$\bar\pi^1\vert_{S_i}\colon S_i\to M$ must be onto on $M$.
\item
There are two possible options for the final step of this algorithm, namely:
\ben
\item
The algorithm ends by giving a submanifold
$S_f\hookrightarrow J^1E$, with $\dim\, S_f\geq m$,
(where \dst S_f=\cap_{i\geq 1}S_i$)
and Euler-Lagrange multivector fields $X_{\Lag}\in\vf^m(S_f)$.
$S_f$ is then called the {\sl final constraint submanifold}.
\item
The algorithm ends by giving a submanifold $S_f$ with $\dim\, S_f<m$,
or the empty set. Then there is no Euler-Lagrange multivector fields
$X_{\Lag}\in\vf^m(S_f)$. In this case the Lagrangian system has no solution.
\een
\eit

{\bf Comment}:
\bit
\item
The problem here considered
can also be treated in an alternative manner by splitting the procedure
into two steps (see \cite{LMM-95} for a more detailed
discussion on this method, using Ehresmann connections):
\ben
\item
First, searching for a submanifold $P_f$ of $J^1E$ and locally decomposable
multivector fields $X_{\Lag}\in\vf^m(J^1E)$, such that:
\ben
\item
$(\inn(X_{\Lag})\Omega_{\Lag})\vert_{P_f}=0$.
\item
$(\inn(X_{\Lag})(\bar\pi^{1*}\omega))\vert_{P_f}=1$.
\item
$X_{\Lag}\colon P_f\to\Lambda^m\Tan P_f$.
\een
This procedure is called the {\sl constraint algorithm}.
\item
Second, searching for a submanifold $S_f$ of $P_f$
and locally decomposable multivector fields
$X_{\Lag}\in\vf^m(J^1E)$, such that:
\ben
\item
$(\inn(X_{\Lag})\Omega_{\Lag})\vert_{S_f}=0$.
\item
$X_{\Lag}$ is a SOPDE on $S_f$.
\item
$X_{\Lag}\colon S_f\to\Lambda^m\Tan S_f$.
\een
\een
It is then clear that, in some particular cases,
locally decomposable and $\bar\pi^1$-transverse
multivector fields, solutions of the evolution equations on some submanifold
$P_f$, can exist, but none of them being SOPDE at any point belonging to $P_f$.
\eit

The local treatment of the singular case shows significative differences to the regular one.
In fact, although the starting equations are the same in both cases
(equations (\ref{eqsG1}) and (\ref{eqsG2})); since the matrix
\dst\left(\frac{\partial^2\lag}{\partial v^A_\nu\partial v^B_\mu}\right)\)
is not regular,  the first group does not imply the SOPDE condition.
Nevertheless, we can impose this condition by making $F^A_\mu=v^A_\mu$,
for every $A,\mu$. Therefore, we have the system of equations
for the coefficients $G^A_{\mu\nu}$ which is (\ref{eqsG}) again.
As we have stated, this system is not compatible in general,
and $S_1$ is the closed submanifold where the system is compatible.
Then, there are Euler-Lagrange multivector fields on $S_1$,
but the number of arbitrary functions on which they depend
is not the same as in the regular case, since it depends on the dimension of
$S_1$ and the rank of the Hessian matrix of $\lag$.
Now the tangency condition must be analyzed in the same way
as in the algorithm of section \ref{ia}

Finally, the question of integrability must be considered.
We make the same remarks as at the end of the above section,
but for the submanifold $S_f$ instead of $J^1E$.

\subsection{Noether's theorem for multivector fields}
\protect\label{ntmvf}

Next we recover the idea of {\sl conserved quantity}
and state Noether's theorem in terms of multivector fields.
First we need a previous lemma:

\begin{lem}
Let $Y\in\vf^m(J^1E)$ be a holonomic multivector field.
If $\zeta\in\df^m(J^1E)$ belongs to the ideal
generated by the contact module, (this ideal is denoted ${\cal I}({\cal M}_c)$), then
$\inn (Y)\zeta =0$.
\end{lem}
\proof
Consider $\bar y\in J^1E$, with $\bar\pi^1(\bar y)=x$,
and a section $\phi\colon M\to E$
such that $j^1\phi$ is an integral section of $Y$ passing
through $\bar y$; that is $j^1\phi (x)=\bar y$.
Then, as $Y$ is SOPDE, there exists $f\in\Cinfty (J^1E)$ with
$$
(\Lambda^m\Tan_xj^1\phi)
\left(\derpar{}{x^1}\wedge\ldots\wedge\derpar{}{x^m}\right) =
(fY)_{\bar y}
$$
and therefore
\beann
(\inn (Y)\zeta )(\bar y) &=& 
(\inn (Y)\zeta )(j^1\phi (x))=
\inn (Y_{j^1\phi (x)})(\zeta (j^1\phi (x)))
\\ &=&
\frac{1}{f(\bar y)} \inn \left( (\Lambda^m\Tan_xj^1\phi)
\left(\derpar{}{x^1}\wedge\ldots\wedge\derpar{}{x^m}\right)\right) (\zeta (j^1\phi (x)))
\\ &=&
\frac{1}{f(\bar y)} \inn \left( \derpar{}{x^1}\wedge\ldots\wedge\derpar{}{x^m}\right)
((j^1\phi)^*\zeta)(x))=0
\eeann
because $\zeta$ belongs to the ideal generated by the contact module and
then $(j^1\phi )^*\zeta =0$.
\qed

Now we are able to prove that:

\begin{teor} {\rm (Noether)}:
Let $\ls$ be a Lagrangian system and $X_{\Lag}\in\vf^m(J^1E)$
an integrable Euler-Lagrange multivector field.
Let $X\in\vf (J^1E)$ be a vector field satisfying the following conditions:
\ben
\item
$X$ preserves the contact module:
$\Lie (X){\cal M}_c\subset{\cal M}_c$.
\item
There exist $\xi\in\df^{m-1}(J^1E)$, and
$\alpha\in{\cal I}({\cal M}_c)\cap\df^m(J^1E)$
such that $\Lie (X)\Lag =\d\xi +\alpha$.
\een
Then $\inn (X_{\Lag})\d (\xi -\inn (X)\Theta_{\Lag}) =0$.
\end{teor}
\proof
First we have that
$$
\Lie (X)\Theta_{\Lag} =\inn (X)\d\Theta_{\Lag} +\d\inn (X)\Theta_{\Lag}
$$
however, as $\Theta_{\Lag}:=\inn({\cal V})\Lag+\Lag\equiv\vartheta_{\Lag}+\Lag$,
$$
\Lie (X)\Theta_{\Lag} =\Lie (X)(\vartheta_{\Lag} +\Lag )=
\Lie (X)\vartheta_{\Lag} +\d\xi +\alpha
$$
and $\Lie (X)\vartheta_{\Lag}\in{\cal I}({\cal M}_c)$,
because $\vartheta_{\Lag}\in{\cal I}({\cal M}_c)$. So we obtain
$$
\inn (X)\d\Theta_{\Lag}+\d\inn(X)\Theta_{\Lag}=\d\xi+\zeta
$$
with $\zeta=\alpha+\Lie(X)\vartheta_{\Lag}\in{\cal I}({\cal M}_c)$. Then
$$
\d (\xi -\inn (X)\Theta_{\Lag})=\inn (X)\d\Theta_{\Lag}-\zeta
$$
and therefore
$$
\inn (X_{\Lag})\d (\xi -\inn (X)\Theta_{\Lag})=
\inn (X_{\Lag})\inn (X)\d\Theta_{\Lag}+\inn (X_{\Lag})\zeta=0
$$
since $\inn (X_{\Lag})\inn (X)\Omega_{\Lag}=0$,
for every $X\in\vf (J^1E)$, because $X_{\Lag}$ is a solution of the
Euler-Lagrange equations, and
$\inn (X_{\Lag})\zeta =0$, by the previous lemma.
\qed

\begin{definition}
Let $\ls$ be a Lagrangian system.
An {\rm (infinitesimal) generalized symmetry} of the Lagrangian system is a
vector field $X\in\vf (J^1E)$ such that it
verifies conditions 1 and 2 of the above theorem.

In particular, $X$ is said to be an {\rm infinitesimal natural symmetry} if
it is the canonical prolongation of a vector field $Z\in\vf (E)$.
\end{definition}

{\bf Remark}:
As a special case, if $\Lie (X)\Lag=0$, then
$\inn (X_{\Lag})\d\inn(X)\Theta_{\Lag}=0$.
These are the so-called {\sl symmetries of the Lagrangian}.

An immediate consequence of the above theorem is the following fact:
if $\alpha=\xi-\inn(X)\Theta_{\Lag}\in\df^{m-1}(J^1E)$,
then $\inn(X_{\Lag})\alpha=0$, and hence,
if $\j\colon N\hookrightarrow J^1E$ is an integral submanifold of $X_{\Lag}$,
we have $\j^*\d\alpha=0$ and then $\d\j^*\alpha=0$.
But the integral submanifolds of $X_{\Lag}$ are images of
canonical lifting of sections, $j^1\phi\colon M\to J^1E$;
therefore
$0=(j^1\phi)^*\d\alpha=\d (j^1\phi)^*\alpha$;
that is:

\begin{teor}{\rm (Noether):}
Let $X$ be an infinitesimal generalized symmetry of the  system
$\ls$. Then, the $(m-1)$-form $\xi-\inn (X)\Theta_{\Lag}$
is closed on the critical sections of the variational problem
posed by $\Lag$.
\end{teor}

\subsection{Example}
\protect\label{example}

Most of the (quadratic) Lagrangian systems in field theories
can be modeled as follows:
$\pi\colon E\to M$ is a trivial bundle (usually $E=M\times\Real^N$)
and then $\pi^1\colon J^1E\to E$ is a vector bundle.
$g$ is a metric in this vector bundle,
${\mit\Gamma}$ is a connection of the projection $\pi^1$,
and $f\in\Cinfty(E)$ is a potential function.
Then the lagrangian function is
$$
\lag (\bar y) =
\frac{1}{2}g(\bar y-{\mit\Gamma}(\pi^1(\bar y)),\bar y-{\mit\Gamma}(\pi^1(\bar y)))+
{\pi^1}^*f(\bar y)
$$
and in natural coordinates takes the form \cite{Sd-95}
$$
\lag =\frac{1}{2}a_{AB}^{\mu\nu}(y)
(v^A_\mu-{\mit\Gamma}^A_\mu(x))(v^B_\nu-{\mit\Gamma}^B_\nu(x))+f(y)
$$
In order to apply the formalism which has been developed in this work,
we consider a model where the matrix of the coefficients $a_{AB}^{\mu\nu}$
is regular and symmetric at every point
(that is, $a^{\mu\nu}_{AB}(y)=a^{\nu\mu}_{BA}(y)$).
This fact is equivalent to the non-degeneracy of the metric $g$.

The associated Poincar\'e-Cartan $(m+1)$-form is then
\beann
\Omega_{\Lag}&=&
-a_{BA}^{\nu\mu}\d v^B_\nu\wedge\d y^A\wedge\d^{m-1}x_\mu-
\\ & &
\derpar{a_{CA}^{\nu\mu}}{y^B}(v^C_\nu-{\mit\Gamma}^C_\nu)
\d y^B\wedge\d y^A\wedge\d^{m-1}x_\mu+
a_{BA}^{\nu\mu}v^A_\mu\d v^B_\nu\wedge\d^mx  +
\\ & &
\left(\derpar{a_{CA}^{\nu\mu}}{y^B}(v^C_\nu-{\mit\Gamma}^C_\nu)v^A_\mu-
\derpar{a_{CA}^{\nu\mu}}{y^B}(v^C_\nu-{\mit\Gamma}^C_\nu)-
\derpar{f}{y^B}  -  a_{AB}^{\mu\nu}\derpar{{\mit\Gamma}_\nu^A}{x^\mu}\right)
\d y^B\wedge\d^mx
\eeann
and the system is obviously regular. Then,
taking (\ref{locsopde}) as the local expression for representatives
of the corresponding classes of Euler-Lagrange multivector fields
$\{ X_{\Lag}\}\subset\vf^m(J^1E)$,
their components $G_{\mu\nu}$ are related by the equations (\ref{eqsG}),
which in this case are
\beq
a_{BA}^{\nu\mu}G^B_{\nu\mu}=
\frac{1}{2}\derpar{a_{CB}^{\nu\mu}}{y^A}
(v^C_\nu-{\mit\Gamma}^C_\nu)(v^B_\mu-{\mit\Gamma}^B_\mu)+\derpar{f}{y^A}+
a_{BA}^{\mu\nu}\derpar{{\mit\Gamma}^B_\nu}{x^\mu}-
\derpar{a_{CA}^{\nu\mu}}{y^B}v^B_\mu(v^C_\nu-{\mit\Gamma}^C_\nu)
\label{eqsGexamp}
\eeq
This system allows us to isolate $N$ of these components
as functions of the remaining $N(m^2-1)$;
and then it determines a family of (classes of) Euler-Lagrange multivector fields.
In order to obtain an integrable class, the condition of integrability
${\cal R}=0$ (where ${\cal R}$ is the curvature of the associated connection)
must hold; that is, equations (\ref{curvcero}) must be added
to the last system.

A simpler case is obtained when ${\mit\Gamma}$ is an integrable connection.
Then there exist local natural charts in $J^1E$ such that ${\mit\Gamma}^A_\mu=0$
\cite{EM-92}; and then
\dst\lag =\frac{1}{2}a_{AB}^{\mu\nu}(y)v^A_\mu v^B_\nu+f(y)\) .
If, in addition, we consider that
the matrix of coefficients is $a^{\mu\nu}_{AB}(y)=\delta_{AB}\delta^{\mu\nu}$
(that is, we take an orthonormal frame for the metric $g$), then we have that
$$
\lag =\frac{1}{2}\delta_{AB}\delta^{\mu\nu}v^A_\mu v^B_\nu+f(y)
$$
In this case equations (\ref{eqsGexamp}) reduce to
$$
\delta_{BA}\delta^{\nu\mu}G^B_{\nu\mu}=\derpar{f}{y^A}
$$
 From this system we can isolate $N$ of the coefficients $G^A_{\mu\nu}$;
for instance, if $\mu,\nu=0,1,\ldots ,m-1$, those for which $\mu=\nu=0$: Thus
$$
\delta_{AB}G^B_{00}=\derpar{f}{y^A}-\sum_{\mu=1}^{m-1}\delta_{AB}G^B_{\mu\mu}
$$
Therefore the Euler-Lagrange multivector fields are
$$
X_{\Lag}=
\bigwedge_{\mu=0}^{m-1}\left(\derpar{}{x^\mu}+v_\mu^A\derpar{}{y^A}+
\delta_\mu^0\left[\derpar{f}{y^A}-
\sum_{\nu=1}^{m-1}\delta_{AB}G^B_{\nu\nu}\right]\derpar{}{v_0^A}+
\sum_{\mu=\gamma\not=0}G_{\mu\gamma}^B\derpar{}{v_\gamma^B}+
\sum_{\mu\not=\gamma}G_{\mu\gamma}^C\derpar{}{v_\gamma^C}\right)
$$
Now, if we look for integrable Euler-Lagrange multivector fields,
the integrability conditions (\ref{curvcero}) must be imposed.

As a particular case, we can consider this last model with
$\dim\,  M=2$. The equations (\ref{eqsG}) are
$$
\delta^{BA}\derpar{f}{y^A}-G^B_{00}-G^B_{11}=0
$$
and the expression of the Euler-Lagrange multivector fields will then be
\beann
X_\Lag &=&
\left(\derpar{}{x^0}+v^A_0\derpar{}{y^A}+G^A_{00}\derpar{}{v^A_0}+
G^A_{01}\derpar{}{v^A_1}\right)
\wedge
\\ & &
\wedge
\left(\derpar{}{x^1}+v^A_1\derpar{}{y^A}+G^A_{10}\derpar{}{v^A_0}+
\left[ G^A_{00}-\delta^{AB}\derpar{f}{y^B}\right]\derpar{}{v^A_1}\right)
\eeann
The integrability conditions (\ref{curvcero}) are
\beann
0 &=& G^B_{01}-G^B_{10}
\\
0 &=&
\derpar{G^B_{10}}{x^0}+v^A_0\derpar{G^B_{10}}{y^A}+G^A_{00}\derpar{G^B_{10}}{v^A_0}+
G^A_{01}\derpar{G^B_{10}}{v^A_1}
\\ & &
-\derpar{G^B_{00}}{x^1}-v^A_1\derpar{G^B_{00}}{y^A}-G^A_{10}\derpar{G^B_{00}}{v^A_0}-
\left( G^A_{00}-\delta^{AC}\derpar{f}{y^C}\right)\derpar{G^B_{00}}{v^A_1}
\\
0 &=&
\derpar{G^B_{00}}{x^0}+
v^A_0\derpar{G^B_{00}}{y^A}-v^A_0\frac{\partial^2f}{\partial {y^A}^2}+
G^A_{00}\derpar{G^B_{00}}{v^A_0}+G^A_{01}\derpar{G^B_{00}}{v^A_1}
\\ & &
-\derpar{G^B_{01}}{x^1}-v^A_1\derpar{G^B_{01}}{y^A}-G^A_{10}\derpar{G^B_{01}}{v^A_0}-
\left( G^A_{00}-\delta^{AC}\derpar{f}{y^C}\right)\derpar{G^B_{01}}{v^A_1}
\eeann
The first relation imposes $G^B_{01}=G^B_{10}$, and then
the others are a system of  $2N$ partial differential equations on the
functions $G^B_{00}$ and $G^B_{01}$. For this system the
hypothesis of Cauchy-Kowalewska's theorem holds \cite{Di-74},
and hence there exist integrable multivector fields solutions of the
Euler-Lagrange equations for $\Lag$.

\section{Conclusions and outlook}

We have studied the integrability of multivector fields
in a differentiable manifold and
the relation between integrable jet fields and multivector fields
in jet bundles, using this to give alternative geometric formulations
of the Lagrangian formalism of classical field theories
(of first order). In particular:
\bit
\item
We have introduced the notion of {\sl integrable multivector fields}
in a manifold $E$, first proving that every orientable
distribution in $E$ is equivalent to a class of
{\sl non-vanishing, locally decomposable} multivector fields,
and then demanding that this distribution be integrable.
The idea of {\sl multiflow} for a certain kind of integrable multivector fields
(which we call {\sl dynamical}) is studied, showing in addition that from
every class of integrable multivector fields we can
select a representative which is dynamical.

In addition, an algorithmic procedure for finding the (maximal) submanifold of
$E$ where a locally decomposable multivector field is integrable
(if it exists) is also developed.

The results obtained can be summarized in the following table:
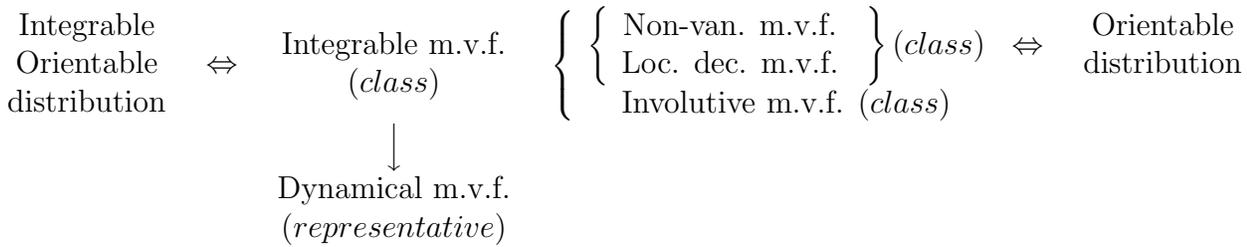
\eit
$$
\begin{array}{cccc}
\begin{array}{c}
\mbox{Integrable}\\
\mbox{Orientable}\\
\mbox{distribution}
\end{array}
&\Leftrightarrow&
\begin{array}{c}
\mbox{Integrable m.v.f.}\\
(class)
\end{array}
&
\left\{
\begin{array}{ccc}
\left\{
\begin{array}{c}
\mbox{Non-van. m.v.f. }\\
\mbox{Loc. dec. m.v.f. }
\end{array}
\right\}
(class)
&\Leftrightarrow&
\begin{array}{c}
\mbox{Orientable}\\
\mbox{distribution}
\end{array}
\\
\mbox{Involutive m.v.f. }(class)
&
\end{array}
\right.
\\
& & \Big\downarrow &
\\
& & \begin{array}{c}
\mbox{Dynamical m.v.f.}\\
(representative)
\end{array}
&
\end{array}
$$
\bit
\item
We have established the relation between integrable jet fields in $J^1E$
and a certain kind of integrable $m$-multivector fields in $E$, by proving that
every orientable jet field $\Psi$ is equivalent to
a class of non-vanishing $m$-multivector fields $\{ Y\}$, and conversely,
in such a way that their associated distributions
${\cal D}(\Psi)$ and ${\cal D}(Y)$ are the same.
Of course, in order to make this equivalence consistent, it is necessary
for the distribution ${\cal D}(Y)$ to be transverse to the 
projection $\pi\colon E\to M$.
So we have:
$$
\left.
\begin{array}{c}
\mbox{Transverse m.v.f.} \\
\mbox{Locally decomp. m.v.f.}
\end{array}
\right\}
(class)
\ \Longleftrightarrow\ 
\mbox{Orientable j.f.}
$$
\item
We have applied these ideas to prove that
every orientable jet field ${\cal Y}$ in $J^1J^1E$ is equivalent to
a class of locally decomposable and non-vanishing $m$-multivector fields
$\{ X\}$ in $J^1E$, and conversely.
In this case, the consistency of this equivalence requires
the distribution ${\cal D}(X)$ to be transverse to the 
projection $\bar\pi^1\colon J^1E\to M$.

The integrability and the SOPDE condition
for multivector fields in this framework are then discussed,
giving several equivalent characterizations for SOPDE
multivector fields, which lead finally to establishing
the equivalence between classes of
holonomic multivector fields and orientable, holonomic jet fields.
So we have the following summarizing scheme:
\eit
$$
\begin{array}{c}
\mbox{Holonomic m.v.f}\\
(class)
\end{array}
\left\{
\begin{array}{ccc}
\left\{
\begin{array}{c}
\mbox{$\bar\pi^1$-Transverse m.v.f.}\\
\mbox{Involutive m.v.f.}
\end{array}
\right\}
(class)
&\Leftrightarrow&
\left\{
\begin{array}{c}
\mbox{Orientable j.f.}
\\
\mbox{Integrable j.f.}
\end{array}
\right.
\\
\left\{
\begin{array}{c}
\mbox{Loc. decom. m.v.f.}\\
\mbox{SOPDE m.v.f.}
\end{array}
\right\}
(class)
&\Leftrightarrow&
\left\{
\begin{array}{c}
\mbox{SOPDE j.f.}
\\
\mbox{Orientable j.f.}
\end{array}
\right.
\end{array}
\right.
\Biggr\}
\begin{array}{c}
\mbox{Holonomic}\\
\mbox{j.f.}
\end{array}
$$
\begin{itemize}
\item
We prove that the evolution equations for first order classical field theories
in the Lagrangian formalism ({\sl Euler-Lagrange equations}\/)
can be written in three equivalent geometric ways:
using multivector fields in $J^1E$, jet fields in $J^1J^1E$
or their associated Ehresmann connections in $J^1E$.
These descriptions allow us to write the evolution equations for field
theories in an analogous way to the dynamical equations for
(time-dependent) mechanical systems.
\item
Using the formalism with multivector fields, we show that
the evolution equations can have no integrable solutions in $J^1E$,
for neither regular nor singular Lagrangian systems.

In the regular case, Euler-Lagrange multivector fields
(that is, SOPDE and solution of the Lagrangian evolution equations)
always exist; but they are not necessarily integrable.
In the singular case, not even the existence of
such an Euler-Lagrange multivector field is assured.
In both cases, the multivector field solution (if it exists) is not unique.
These features are significant differences in relation to the analogous situation in mechanics.
\item
With regard to these facts, for the singular case we outline an algorithmic procedure
in order to obtain the submanifold of $J^1E$ where
Euler-Lagrange multivector fields exist.

Further research on these latter topics will be carried out in the future.
\item
A version of Noether's theorem using multivector fields is also proved.
\eit

In a forthcoming paper devoted to the Hamiltonian formalism
of field theories, we will extend these results to the dual bundle
of $J^1E$.

\appendix

\section{Canonical structures of a jet bundle}

(Following \cite{Gc-73}. See also \cite{EMR-96}).

Consider a first-order jet bundle
$J^1E\mapping{\pi^1}E\mapping{\pi}M$.
Let $\phi\colon M\to E$ be a section of $\pi$, $x\in M$ and $y=\phi(x)$.
The {\sl vertical differential} of the section $\phi$ at the point $y\in E$
is the map
$$
\begin{array}{ccccc}
\d^v_y\phi&\colon&\Tan_y E & \longrightarrow & {\rm V}_y(\pi) \\
  & & u & \longmapsto & u-\Tan_y(\phi\circ\pi)u
\end{array}
$$
If $(x^\mu ,y^A)$ is a natural local system of 
$E$ and $\phi=(x^\mu ,\phi^A(x^\mu ))$, then
$$
\d^v_y \phi \left(\derpar{}{x^\mu}\right)_y =
-\left(\derpar{\phi^A}{x^\mu}\right)_y\left(\derpar{}{y^A}\right)_y
\quad , \quad
\d^v_y \phi \left(\derpar{}{y^A}\right)_y =
\left(\derpar{}{y^A}\right)_y
$$
Observe that the vertical differential splits $\Tan_y E$
into a vertical component and another which is tangent
to the imagen of $\phi$ at the point $y$. 
In addition, $\d^v_y \phi (u)$ is the projection of $u$
on the vertical part of this splitting.

\begin{definition}
Consider $\bar y\in J^1E$ with
$\bar y\stackrel{\pi^1}{\mapsto}y\stackrel{\pi}{\mapsto}x$
and $\bar u\in\Tan_{\bar y}J^1E$.
The {\rm structure canonical form} of $J^1E$ is a $1$-form $\theta$
in $J^1E$ with values on $\pi^{1^*}{\rm V}(\pi)$ which is defined by
$$
\theta (\bar y;\bar u):=(\d^v_y \phi)(\Tan_{\bar y}\pi^1 (\bar u))=
(\Tan_{\bar y}\pi^1 -\Tan_{\bar y}(\phi\circ\bar\pi^1))(\bar u)
$$
where the section $\phi$ is a representative of $\bar y$.
\label{theta}
\end{definition}

\begin{definition}
Let $\theta$ be the structure canonical form of $J^1E$.
It can be considered as a $\Cinfty (J^1E)$-linear map
$\theta\colon\Gamma (J^1E,\pi^{1^*}{\rm V}(\pi))^*\to\df^1(J^1E)$.
The image of this map is called the {\rm contact module}
or {\rm Cartan distribution} of $J^1E$.
It is denoted by ${\cal M}_c$.
\end{definition}

The expression of $\theta$ in a natural local system of $J^1E$ is
\dst\theta=({\rm d}y^A - v^A_{\mu}{\rm d}x^{\mu})\otimes \derpar{}{y^A}\) .
Then ${\cal M}_c$ is generated by the forms
$\theta^A=\d y^A-v^A_\mu\d x^\mu$ in this open set.

If $\phi\colon M\to E$ is a section of $\pi$,
we denote by $j^1\phi\colon M\to J^1E$ its canonical prolongation.
Then, a section $\psi\colon M\to J^1E$ is said to be {\sl holonomic}
iff $\psi=j^1\phi$, for some $\phi$.
Holonomic sections can be characterized using the
structure canonical form or, equivalently,
the contact module, as follows:

\begin{prop}
Let $\psi\colon M \to J^1E$ be a section of $\bar\pi^1$.
The following assertions are equivalent:
\ben
\item
$\psi$ is a holonomic section.
\item
$\psi^*\theta=0$.
\item
$\psi^*\alpha =0$, for every $\alpha\in {\cal M}_c$.
\een
\label{holsec}
\end{prop}

\begin{definition}
\ben
\item
For every $\bar y\in J^1E$, consider the canonical isomorphism
$$
{\cal S}_{\bar y}\colon
\Tan_{\bar\pi^1 (\bar y)}^*M\otimes{\rm V}_{\pi^1 (\bar y)}(\pi)
\longrightarrow {\rm V}_{\bar y}(\pi^1 )
$$
which consists in associating with an element
$\alpha\otimes v\in
\Tan_{\bar\pi^1 (\bar y)}^*M\otimes{\rm V}_{\pi^1 (\bar y)}(\pi)$
the directional derivative in $\bar y$ with respect to $\alpha\otimes v$.
Taking into account that  $\alpha\otimes v$ acts in $J^1_yE$
by translation
$$
{\cal S}_{\bar y}(\alpha\otimes v ):= D_{\alpha\otimes v}(\bar y)
\colon f\mapsto \lim_{t\to 0}\frac{f(\bar y+t(\alpha\otimes v))-
f(\bar y)}{t}
$$
for $f\in\Cinfty (J^1_yE)$.
Then we have the following isomorphism of $\Cinfty (J^1E)$-modules
$$
{\cal S}\colon
\Gamma (J^1E,\bar\pi^{1^*}\Tan^*M)\otimes\Gamma (J^1E,\pi^{1^*}{\rm V}(\pi))
\longrightarrow \Gamma (J^1E,{\rm V}(\pi^1 ))
$$
which is called the {\rm vertical endomorphism} ${\cal S}$.
\item
As ${\cal S}\in\Gamma (J^1E,(\pi^{1^*}{\rm V}(\pi))^*)\otimes
\Gamma (J^1E,{\rm V}(\pi^1 ))\otimes\Gamma(J^1E,\bar\pi^{1^*}\Tan M)$
(where all the tensor products are on $\Cinfty (J^1E)$);
the {\rm vertical endomorphism} ${\cal V}$ arises from
the natural contraction between
the factors $\Gamma (J^1E,(\pi^{1^*}{\rm V}(\pi))^*)$ of ${\cal S}$
and $\Gamma (J^1E,\pi^{1^*}{\rm V}(\pi))$ of $\theta$; that is
$$
{\cal V}=
\inn ({\cal S})\theta\in\df^1(J^1E)\otimes\Gamma (J^1E,{\rm V}(\pi^1 ))
\otimes\Gamma (J^1E,\bar\pi^{1^*}\Tan M)
$$
\een
\end{definition}

If $(x^\mu ,y^A,v^A_\mu )$ is a natural system of coordinates,
the local expression of ${\cal S}$ and ${\cal V}$ are given in this local system by
$$
{\cal S}=\xi^A\otimes\derpar{}{v^A_\mu}\otimes\derpar{}{x^\mu}
\quad ;\quad
{\cal V}=\left(\d y^A-v^A_\mu\d x^\mu\right)\otimes
\derpar{}{v^A_\nu}\otimes\derpar{}{x^\nu}
$$
where $\{\xi^A\}$ is the local basis of 
$\Gamma (J^1E,\pi^{1^*}{\rm V}(\pi))^*$
which is dual of \dst\left\{\derpar{}{y^A}\right\}\) .

\subsection*{Acknowledgments}

We wish to thank Mr. Jeff Palmer for his assistance in
preparing the English version of the manuscript.

We are also grateful for the financial support
of the CICYT TAP97-0969-C03-01.

\end{document}